\documentclass[preprint,11pt]{elsarticle}

\usepackage{amssymb}

\usepackage[ansinew,latin1]{inputenc}
\usepackage[english]{babel}
\usepackage{amsfonts} 
\usepackage{amsmath}
\usepackage{amssymb}
\usepackage{latexsym}
\usepackage{amsthm}
\usepackage{verbatim}      
\usepackage{bm,enumerate}
\usepackage[normalem]{ulem}     
\usepackage{cancel}    
\usepackage{lineno}     
\usepackage{multirow}
\usepackage{array}
\usepackage{color}
\usepackage{fullpage}
  
\addto\captionsenglish{}

\def\mass{\mathbb{M}}
\def\cgmass{\overline{\mathbb{M}}}

\def\PUSHFWD#1{\COP_*{#1}}                 
\def\PUSHFWDA#1{\COP^{(1)}_*{#1}} 
\def\PUSHFWDB#1{\COP^{(2)}_*{#1}} 
\def\COPINV{\COP^\dagger}                  
\def\PULLBCK#1{\COPINV_*{#1}} 
\def\TI{\mathrm{I}}
\def\TII{\mathrm{II}}
\def\TIII{\mathrm{III}}

\def\GCOP{\COP^\sharp} 

\def\Bb{\mathcal{B}}

\def\Hh{\mathcal{H}}

\def\Ll{\mathcal{L}}

\def\Nn{{\mathcal{N}}}

\def\Rr{{\mathcal{R}}}

\def\E{\mathbb{E}}

\def\R{\mathbb{R}}

\def\EXP#1{e^{#1}}

\def\EXPECT{{\mathbb{E}}}

\def\RELENT#1#2{\Rr\left(#1|#2\right)}

\def\VAR{\mathrm{Var}}

\def\SCPROD#1#2{\langle {#1},{#2}\rangle}

\def\COMMA{\,,}             
\def\PERIOD{\,.}            
\def\SEP{{\,|\,}}           

\def\VIZ#1{(\ref{#1})}      

\def\BIGO{\mathcal{O}}

\def\DIAG{\mathrm{diag}\,}

\newtheorem{theorem}{Theorem}[section]

\newtheorem{corollary}{Corollary}[section]

\newtheorem{remark}{Remark}[section]

\theoremstyle{definition}
\newtheorem*{prf}{Proof}

\newtheorem{thm}{Theorem}

\def\PrS{\COP \Sigma(x)\COP^{tr}}

\def\bx{\BARIT{x}}
\def\bp{\BARIT{p}}
\def\bq{\BARIT{q}}
\def\tq{\tilde{q}}
\def\tp{\tilde{p}}
\def\tb{\tilde{b}}
\def\bbr{\BARIT{b}}
\def\BF{\BARIT{F}}
\def\BP{\BARIT{P}}
\def\bX{\BARIT{X}}
\def\bB{\BARIT{B}}
\def\tX{\tilde{X}}

\def\IP{U}

\def\BDelta{\BARIT{\Delta}}

\def\ba{\mathbf{a}}
\def\BPhi{\mathbf{\Phi}}

\def\BARIT#1{{\bar {#1}}}

\def\COP{\mathbf{\Pi}}

\def\PXT{\{X_t\}_{t\geq 0}}
\def\tPXT{\{\tilde X_t\}_{t\geq 0}}
\def\CGPXT{\{\BARIT X_t\}_{t\geq 0}}
 \def\PATHPT{P_{[0,T]}}
 \def\tPATHPT{Q^{\theta}_{[0,T]}}
 \def\CGPATHPT{\BARIT Q^{\theta}_{[0,T]}}
\def\PATHP{P}
\def\tPATHP{Q^{\theta}}
\def\tPATHPOP{Q^{\theta^*}}
\def\CGPATHP{\BARIT Q^{\theta}}

\def\0o{\mathbb{O}}

\def\LAWEXACT{P}
\def\LAWAPPROX{Q}

\def\PATHS{\LAWEXACT_{[0,T]}}
\def\PATHSAPP{\LAWAPPROX_{[0,T]}}

\def\RELENT#1#2{\mathcal{R}\left({#1}\SEP{#2}\right)}
\def\ENTRATE#1#2{\mathcal{H}({#1}\SEP{#2})}

\def\EXPECT{{\mathbb{E}}}
\def\EXPECTWRT#1#2{\E_{#1}\left[{#2}\right]}

\begin{document}
\begin{frontmatter}
%
%
\title{Path-space variational inference for non-equilibrium coarse-grained systems}

\author[uoc,forth]{Vagelis Harmandaris}
\ead{harman@uoc.gr}
\author[uoc]{Evangelia Kalligiannaki\corref{co1}}
\ead{ekalligian@tem.uoc.gr}
\author[umass]{Markos Katsoulakis\corref{co2}}
\ead{markos@math.umass.edu}
\author[udel]{Petr Plech\'a\v{c}} 
\ead{plechac@math.udel.edu}
\address[udel]{Department of Mathematical Sciences,
University of Delaware,
Newark, Delaware, United States }
\address[umass]{Department of Mathematics and Statistics, University
of Massachusetts at Amherst, United States}
\address[uoc]{ Department of Mathematics and Applied Mathematics, 
University of Crete, Greece}
\address[forth]{Institute of Applied and Computational Mathematics (IACM), Foundation for Research and Technology Hellas
	(FORTH), IACM/FORTH, GR-71110 Heraklion, Greece}

\cortext[co1]{Corresponding author}
\cortext[co2]{Principal corresponding author}
%
\begin{abstract}
In this paper we discuss information-theoretic tools for obtaining optimized coarse-grained molecular 
models for both equilibrium and non-equilibrium molecular simulations. 
The latter are ubiquitous in physicochemical and biological applications, where they are typically associated with coupling mechanisms, multi-physics and/or boundary conditions.
In general the non-equilibrium steady states are not known explicitly as they do not necessarily 
have a Gibbs structure. 

The presented approach can compare microscopic behavior of molecular systems to parametric and
non-parametric coarse-grained models  using the relative entropy between distributions on 
the path space and setting  up a corresponding path-space variational inference  problem. 
The methods can become entirely data-driven when the microscopic dynamics are replaced with 
corresponding correlated data in the form of time series. 
Furthermore, we present connections and generalizations of force matching methods in 
coarse-graining with  path-space information methods. We demonstrate the  enhanced transferability 
of information-based parameterizations to   different observables, at a specific thermodynamic point,    
due to information  inequalities.

We discuss methodological connections  between information-based coarse-graining
 of molecular systems  and variational inference methods primarily developed in the machine learning community. 
However, we note that the work presented here addresses variational inference  for  correlated time series due to the focus on dynamics. The applicability of the proposed methods is demonstrated on  high-dimensional stochastic processes 
given by overdamped and driven Langevin dynamics of interacting particles. 
\end{abstract}


\end{frontmatter}

\newpage
\section{Introduction}\label{Intro}
Molecular dynamics simulations at microscopic (e.g. atomistic) level have 
capability of providing quantitative information about rheological, mechanical, chemical and electrical properties of molecular systems, \cite{Larson1999,Doi-Edwards1986}. 
However, the enormous range of length and time scales involved in such complex materials
presents a challenging computational task, in particular, due to a wide disparity of relaxation times. 
 
A standard methodology  in order to overcome problems of long relaxation times of complex systems is to 
abandon the chemical detail and describe the molecular system by fewer (the most relevant) degrees of freedom.
The choice of the latter depends entirely on the physical problem under question. 
Such particle-based, systematic {\it coarse-grained} (CG) models of molecular systems are developed by averaging out the details at the molecular level, and by representing groups of atoms by a single CG particle. 
Then the effective coarse-grained interaction potentials (more precisely free energies) are derived from the microscopic details of the atomistic model. 
The coarse-grained potentials and force fields can be derived through different methods, such as 
the inverse Boltzmann method, force matching and relative entropy, \cite{LyubLaa2004, g219, Soper1996, Harmandaris2006a, IzVoth2005, IzVoth2005a, Shell2008, Shell2009}.
Applying these methods in the context of best-fit procedures in parametrized families of CG models the structural properties of systems at {\it equilibrium} can be described with accuracy which is related to a metric 
used for the parameter fitting procedure. 
However, the  above mentioned coarse-graining parametrization techniques do not address dynamical 
properties of the model and are restricted to systems already at a (equilibrium) Gibbs state. 
 
Furthermore, there are several important issues related to systematic CG models using microscopic information for molecular systems under non-equilibrium conditions: (a) the whole approach is based on the fact that there is a direct connection between structural properties (like pair distribution functions) and CG interaction potentials; i.e., the renormalization group map or Boltzmann relation, see for instance
\cite{LyubLaa2004, kremer, Harmandaris2006a}. This is certainly true at equilibrium and 
near to equilibrium but may not be the case for systems far from equilibrium; 
(b) since the CG interaction potential intrinsically involves entropy, it is not clear 
what is the dependence of the effective CG force field with respect to the external forces (if they are present); (c) predicting the dynamics (or incorporating the proper friction in the equations of motion) in the CG non-equilibrium model is not clear, \cite{Fritz2011a,Baig2010}. All these aspects are, in principle, relevant in any application of a systematic CG 
model for a molecular system under non-equilibrium conditions.

Recently, several methods for coarse-graining of stochastic models based on information theory have 
 appeared in the literature, ~\cite{shell1, Bilionis:2012, zabaras}.
These methods employ entropy-based techniques that estimate discrepancy between (probability) measures.  
Using entropy-based analytical tools has
proved essential for deriving rigorous results for passage from interacting particle
models to mean-field description, e.g., \cite{Landim}.
Applications of these methods  to the error analysis of coarse-graining of stochastic particle
systems have been introduced in \cite{KV, KPRT1,KPRT2,KPRT3,KMV,KMV2}. 
Independently of such rigorous mathematical work, the engineering community developed
entropy-based computational techniques that are used for constructing approximations of coarse-grained 
potentials for models of large biomolecules and polymeric systems (fluids, melts),
where 
the optimal parametrization of effective potentials is based on minimizing the relative entropy between {\em equilibrium} Gibbs states, e.g., \cite{shell1,shell2, shell3, Bilionis:2012, zabaras}. 
Note, that other works in the literature are primarily based on observable-matching using either structural distribution functions, such as the inverse Boltzmann method, 
inverse Monte Carlo  methods, \cite{Soper1996,LyubLaa2004,g219}, or  
averaged forces on CG particles \cite{IzVoth2005a, IzVoth2005}. 
These methods were used with a great success in coarse-graining of macromolecules, 
see, e.g., \cite{kremer,mp, Harmandaris2006a,Harmandaris2007a,Harmandaris2009a,Fritz2011a,Baig2010}. 
Recent review articles \cite{Fritz2011a, NoidReview:2011,VothReview:2013} give a detailed overview 
of coarse-graining techniques applied to systems at equilibrium. 
Finally, effective coarse equilibrium dynamics for systems with temporal scale separation  modeled by overdamped Langevin dynamics were   studied in \cite{Lelievre2010}.

  Evolution of coarse-grained variables corresponding to Hamiltonian  microscopic  dynamics can be described exactly with  the Mori-Zwanzing formalism  leading  to 
a stochastic integro-differential system  with strong memory terms, known as the generalized Langevin equation  (GLE) \cite{Zwanzig1961, Mori1965}, that is in principle computationally intractable. 
Therefore either a scaling of CG dynamics or approximations of the GLE are used 
\cite{Ottinger, Harmandaris2009a, Harmandaris2009b, Fritz2011a, Briels2001, chorin1, Buscalioni2010, Guenza2011, DarvePNAS2009, Karniadakis2015}.

  Approximate dynamical models  in a parametrized formulation have also been considered  in recent studies, most of them based on the well established equilibrium parametrization methods described above.
For example, authors in \cite{VothDyn2006, VothDyn2015}  propose optimal  CG  parametrized Langevin dynamics based on the force matching method.

 In order to extend the information-theoretic approach developed in \cite{Shell2008} for coarse-graining of Gibbs states to dynamics, a parameter fitting procedure for dynamic coarse-grained models was developed in \cite{EspanolZuniga2011}. The method proposed there is based on minimization of the relative entropy between discrete two-time transition probabilities associated with the diffusion process, in this case Langevin dynamics.  
 Moreover,
authors demonstrate  that the  relative entropy minimization  can be interpreted as a force-matching problem.
 The use of the two-time step probability limits the applicability of the approach  to short time dynamics
while the  discretization time step appears explicitly in the minimization problem leading to time step depended optimal parameters.
  In a recent article,   \cite{Dequidt2015}, authors  attempt to overcome the short time limit using Bayesian inference to identify most probable parameters   for a given time series of microscopic states, i.e. in a path space perspective.
  The authors provide also the connection with force-matching where  though again  the  optimal parameter set depends on the time step of the numerical discretization scheme.  
 The  {\it relative entropy rate} (RER) functional  for Markov Chains proposed in \cite{PK2013} and \cite{KP2013} is similar to the functional that defines the best-fit optimization in \cite{EspanolZuniga2011}, being the relative entropy per unit time for stationary processes.  
  The formulation of path-space  relative entropy for continuous time process in the present work illustrates   that both path-space relative entropy and  RER 
are independent of the time step  for any numerical discretization scheme. This fact is further  demonstrated in \ref{numerics-sec5} where we study the RER minimization for numerical schemes for the Langevin and the overdamped Langevin dynamics. 

In fact, in the present article, as well as earlier in \cite{KP2013}, we show that path-space information metrics 
such as RER provide  a  general framework in several additional  directions: 
(a) they are applicable to discrete-time or discretized (as in \cite{EspanolZuniga2011}) dynamics and  also to broad classes of continuous time stochastic dynamics 
such as Kinetic Monte Carlo algorithms (e.g., reaction-diffusion mechanisms on lattices),
 biochemical reaction networks and semi-Markov processes, \cite{PK2013}, \cite{KP2013}, \cite{PVK2013}, 
(b)  apply to non-equilibrium problems without analytically known Gibbs states and a detail balanced condition 
(irreversible processes),
including driven Langevin models, Kinetic Monte Carlo algorithms and reaction networks, and most importantly,
(c) our  RER perspective shows that the corresponding optimization methods are  extendable to both finite and infinite times, guaranteeing  quantified predictive capability even at long time regimes, \cite{DKPP}.  
In order to demonstrate the abstract principles of the information-theoretic framework on  the {\it path space}
we present coarse-graining of dynamics described as solutions to Ito's stochastic differential equations. 
The presented information-theoretic methodology allows us to  build optimal coarse-grained dynamics  as systematic approximations of microscopic stochastic dynamics in the same general class of Markovian dynamics, e.g.,  of stochastic differential equations.  
The path-space information approach  compares microscopic and coarse-grained  dynamics using the relative entropy between distributions on the path space (see Theorem~\ref{REminThm}) and sets up the corresponding path-space  variational inference (parameter optimization) problems. 

One of the mathematical and computational  novelties of the presented approach lies in the derivation of 
 {\em path space force-matching} conditions  which are applicable to both equilibrium and non-equilibrium systems.
This path-space information theory  formulation provides a natural generalization to  non-equilibrium systems 
of the force matching methods developed earlier for systems at equilibrium, \cite{IzVoth2005, KHKP1}. 
Moreover, we demonstrate here the equivalence of the relative entropy rate  (RER)  and 
force-matching type optimization methods, in analogy  to the  equilibrium case, studied in  
\cite{Voth2008a, NoidReview:2011, KHKP1}.
Furthermore, due to  information inequalities such as \VIZ{CKP} and \VIZ{DKP-bound}, 
we have that path information-based coarse-graining implies transferability of the parameterization to  all reasonable  observables by only training  a single observable,   at the specific thermodynamic point,  namely the relative entropy.
We explore a different but asymptotically equivalent (in the number of data) perspective, 
which is data-driven in the sense  that it treats the microscopic simulator 
only as means of producing statistical data in the form of time series.

We also stress in this work connections between  information-based work for coarse-graining  
and  the variational inference methods primarily developed in the machine learning community. 

In Section~\ref{EqRE} we give a short overview of variational inference  and coarse-graining  
as well as connections to machine learning methods and computational approaches.
Next in Section~\ref{AtCG} we describe the   microscopic   and coarse-grained   models  and the corresponding dynamics. Section~\ref{Sec:VarInf} describes  the path space relative entropy  and relative entropy rate minimization problems and sets up the theoretical tools for comparing the microscopic and coarse dynamics at finite and 
long-time, stationary regimes.
In Section~\ref{RERminS} we prove that the minimization of the  path-space relative entropy is a path space force matching problem providing a generalization of the known equilibrium force matching method to non-equilibrium systems. 
Applications of  our results are presented in Section~\ref{Langevin-sec}  for the  Langevin  dynamics demonstrating applicability of the proposed method to molecular systems while presenting specific examples of coarse graining transformations. Moreover, a  direct study of the relative entropy rate minimization for the 
discretized Langevin and overdamped Langevin dynamics is discussed,  verifying the validity of the continuous time optimal coarse-grained model for the corresponding time discretized  schemes.
In Section~\ref{data_driven} we present the use of the relative entropy minimization as a means of 
optimizing the information content in a coarse model with respect  to available  time series data  coming from a fine-scale simulation. Finally in the last section we summarize the contributions  of the present work.


\section{Variational inference methods and coarse-graining of molecular systems.}\label{EqRE}

Here  we give a short overview of information-based coarse-graining of molecular systems and highlight  connections with variational inference methods primarily developed in the machine learning literature and discuss computational approaches.

\subsection{Information-theoretic methods for coarse-graining molecular systems.}\label{InfoCG}
Computational methods developed for parameterizing coarse-grained models at equilibrium, such as 
inverse Monte Carlo, \cite{LyubLaa2004, g219},  inverse Boltzmann, \cite{kremer, Soper1996,g219},  
force matching,  \cite{IzVoth2005, IzVoth2005a}, and relative entropy, \cite{Shell2008, Shell2009},
provide a development of the CG interaction potential by considering a pre-selected  set of observables 
$\phi_i$, $i=1,...,\ell$ and 
then minimizing a fitting functional over the parameter space $\Theta\subseteq \R^k$.
Typical choices of observables $\phi_i$ are radial distribution functions, \cite{kremer,mp}, and forces between CG particles,
\cite{IzVoth2005a, IzVoth2005,VothReview:2013}.
  A family of parametrization methods  is  described by 
\begin{equation}\label{IMC1}
\min_{\theta \in \Theta} \sum_{i=1}^\ell |\EXPECT_\mu[\phi_i] - \EXPECT_{\bar \mu^{\theta}}[\phi_i]|^2\, ,
\end{equation}
where $\mu$ denotes the fine-scale Gibbs equilibrium distribution and $\bar\mu^{\theta}$ the  parametrized coarse-grained Gibbs distribution. Such methods are referred to as the  iterative   inverse Boltzmann method
and the inverse Monte Carlo methods, \cite{LyubLaa2004, g219,kremer, Soper1996}.  
Clearly any parametrization based on a minimization  principle such as 
 \VIZ{IMC1}  depends on the specific choice of observables, 
while the accurate simulation of other observables, which are not part of the parameterization \VIZ{IMC1} is not necessarily guaranteed. 

Before we continue further we introduce mathematical concepts involved in coarse graining. 
In abstract terms we consider the original (microscopic) model 
defined on a measurable space $(\Omega, \Bb)$, where $\Omega$ represents the state (configuration) space and
$\Bb$ denotes the $\sigma$-algebra on $\Omega$, and the coarse-grained model on $(\bar \Omega, \bar \Bb)$ with
the coarse-graining map 
\begin{equation}
\label{CG-map1}
\COP:\Omega \to \bar\Omega\, .
\end{equation}
The elements of the coarse state space $\bar\Omega$ (the coarse degrees of freedom) are 
thus $\bar\omega = \COP \omega$.
We use the bar ``$\bar{\ }$'' notation for objects related to the coarse-grained model. 
The (probability) measures on
the microscopic space $(\Omega,\Bb)$ are mapped (pushed-forward) by the map $\COP_*$, $\bar P = \COP_* P$,
\begin{equation}\label{defPFW}
\EXPECT_{\bar P}[\bar\phi] \equiv \int_{\bar\Omega} \bar\phi\; d(\PUSHFWD{P}) = \int_{\Omega} \bar\phi\circ \COP\, dP\COMMA
\end{equation}
where $\bar\phi: \bar\Omega\to\R$,
or equivalently $(\PUSHFWD{P})(\bar B) = P(\COP^{-1}(\bar B))$ where $\bar B \in \bar\Bb$.

The relative entropy (Kullback-Leibler divergence),  \cite{Cover}, of two probability measures $P(d\omega)$ and
$Q(d\omega)$ on a common measurable space $(\Omega,\Bb)$ is given by  
\begin{equation}\label{relent1}
\RELENT{P}{Q} = \int_\Omega \log\frac{dP(\omega)}{dQ(\omega)}\,P(d\omega)\COMMA
\end{equation}
provided $P\ll Q$, i.e., $P$ is absolutely continuous with respect to $Q$,
and $\RELENT{P}{Q}=+\infty$ otherwise.
The functional $\RELENT{P}{Q}$ defines a pseudo-distance between two measures 
as $\RELENT{P}{Q}\geq 0$ and $\RELENT{P}{Q}=0$
if and only if $P=Q$, $P$-a.s. In the case these probability measures have
corresponding probability densities  $p(\omega)$ and $q(\omega)$ relation \VIZ{relent1}
 becomes $\RELENT{P}{Q} = \int_\Omega \log\frac{p(\omega)}{q(\omega)}\,p(\omega)d\omega$. 

In contrast to the observable-centered perspective of (\ref{IMC1}), we turn our attention to information inequalities and their implications for coarse-graining. For instance,  the Csisz\'ar-Kullback-Pinsker (CKP) inequality, \cite{Cover}, when applied to an observable 
$\phi$ and the fine-scale and coarse-grained  equilibrium distributions considered in \VIZ{IMC1} readily gives rise to an error bound
\begin{equation}\label{CKP}
|\EXPECT_\mu[\phi] - \EXPECT_{\bar \mu^{ \theta}}[\phi]|
\leq ||\phi||_\infty \sqrt{2 \RELENT{\PUSHFWD{\mu}}{\bar\mu^{\theta}}}\, ,
\end{equation}
where $\COP$ is the coarse-graining map, $\PUSHFWD{\mu}$ denotes the push-forward 
of the microscopic measure $\mu$ defined above,   and  $\| \phi\|_\infty = \sup\{|\phi(x)|: x\in \Omega\} $   is the uniform norm of $\phi$. 

Since the microscopic and coarse-grained problems are formulated on different measure spaces it is necessary to 
clarify how the relative entropy is defined in (\ref{CKP}). 
Note that 
from the computational point of   view  $\RELENT{\PUSHFWD{\mu}}{\bar\mu^{\theta}}$ is still computed on the
microscopic space $\Omega$ using the formula \VIZ{defPFW} for the push-forward. We also refer to Section~\ref{data_driven} for precise formulas and related estimators.
Another option for comparing the two models in the terms of relative entropy is by mapping the coarse model
on the microscopic space, i.e., defining a (microscopic) reconstruction map. 
The reconstruction map can be viewed as a generalized inverse 
$\COPINV:\bar\Omega \to \Omega$, 
and as such is not necessarily unique. Associated with the reconstruction
map is the pull-back map $\PULLBCK{\bar P}$ that defines a measure on the microscopic space $\Omega$.
In that way we can also compare the models by considering  $\RELENT{P}{\PULLBCK{\bar Q^\theta}}$, a strategy we apply in Section~\ref{RERminS}.

Returning to (\ref{CKP}),  as was   pointed out in \cite{KV} in the context of coarse-graining
of stochastic  lattice systems and Kinetic Monte Carlo algorithms, the CKP inequality provides a strong indication that the relative entropy can control all observables $\phi$. 
In our context,  the  CKP inequality  \VIZ{CKP} demonstrates  that minimizing the {\em single observable} given by the relative entropy $\RELENT{\PUSHFWD{\mu}}{\bar\mu^{\theta}}$, 
as  proposed  in \cite{shell1,shell2, zabaras, NoidReview:2011}, i.e., training the parametric coarse-grained model based on
\begin{equation}\label{ENT1}
\min_{\theta\in \Theta} \RELENT{\PUSHFWD{\mu}}{\bar\mu^{\theta}}\COMMA
\end{equation}
instead of \VIZ{IMC1}, will provide reliable coarse-graining parameterizations, applicable to  {\em various observables}.  
Hence, due to  \VIZ{CKP}, information-theoretic methods  give rise to enhanced {\em transferability} properties of the resulting 
coarse-grained model with respect to other observables $\phi$,  t at a specific thermodynamic point.

\begin{remark}
{\rm
Very  recently a sharper version of the CKP inequality has been established in \cite{Dupuis}. 
Furthermore, in \cite{DKPP}, 
it is shown that such information inequalities  can be extended to path-space observables, e.g.,
ergodic averages, correlations, etc. The error between averaging an observable under the probabilistic model
described by $P$ as compared to the model given by $Q$ can be then bounded 
\begin{equation}\label{DKP-bound}
 \Phi_{-}(P, Q; \phi)
\leq \EXPECT_{P}[\phi] - \EXPECT_{Q}[\phi] \leq
 \Phi_+(P, Q; \phi)\PERIOD
\end{equation}
In \cite{DKPP} the authors refer to  $\Phi_{\pm}(P, Q; \phi)$ as ``goal-oriented divergence'' because it has the properties 
of a divergence both in the sense of probabilities $P$ and $Q$  and observables $\phi$:
$\Phi_{+}(P, Q; \phi)\ge 0$, (resp. $\Phi_{-}(P, Q; \phi) \le 0$) and  $\Phi_{\pm}(P, Q; \phi)=0$ if and only if $P=Q$ a.s. 
or $\phi$ is constant $P$-a.s.
Furthermore,  $\Phi_{\pm}(P, Q; \phi)$  admits an expansion in the (small)  relative entropy, \cite{DKPP}:
$$
\Phi_{\pm}(P, Q; \phi)=\pm\sqrt{\VAR_P[\phi]}\sqrt{2\RELENT{P}{Q}} + \BIGO(\RELENT{P}{Q})\COMMA
$$ 
that captures the key properties of this new divergence. 
Similar relations hold for path-space observables, with the role of relative entropy  played by the {\em relative entropy rate} (RER), \cite{DKPP},
defined below in this section.
}
\end{remark}

\subsection{Variational inference and machine learning.}
At this point we digress and point out connections between this and earlier information-based 
works for coarse-graining  
and  variational inference methods. 
In the context of statistical variational inference 
(see for instance, \cite{Wainwright:2008,MacKay:2003,Murphy:2012})
one defines a flexible and rich enough parametrized family of distributions,  
and one finds the member of the family 
which is closest to the posterior distribution in a suitable metric. 
Subsequently one samples from that approximating distribution instead of the posterior itself. 
In this sense the inference problem is tackled using an optimization principle, hence  the term ``variational inference''.  

More precisely, if we denote the  posterior by $P$ and the parametrized family by $Q^\theta$ over the parameter space  
$\theta \in \Theta\subseteq \R^k$ we typically consider in variational inference their ``distance'' as the minimization 
over all  corresponding relative entropies 
\begin{equation}\label{VP1}
\min_{\theta\in \Theta} \RELENT{P}{Q^\theta}\, .
\end{equation}
Reversing the order to $\RELENT{Q^\theta}{P}$ can also be considered (as we also do in this paper), capturing different aspects of the posterior $P$ due to the non-symmetry of the relative entropy, \cite{Murphy:2012,MacKay:2003}. In addition to the single parameter vector  optimization 
in (\ref{VP1}), we can also  consider  approximating Bayesian inference. This latter perspective  gives  rise to variational Bayesian inference, also known as ensemble learning, \cite{MacKay:2003}.

In many important cases, e.g., when $Q^\theta$ is a class of mean field models, that is a parametrized family of product distributions over the state space, the optimization problem  (\ref{VP1}) can be solved analytically, 
\cite{Wainwright:2008,MacKay:2003,Murphy:2012}. 
For instance, the mean field theory  for  the ferromagnetic Ising model can be also derived as 
the solution of (\ref{VP1}), when $Q^\theta$ is a family of product distributions, \cite{MacKay:2003}. 
Furthermore, more complex parametric families $Q^\theta$ than mean-field have been also considered in the literature, see for instance \cite{Wainwright:2008}. 
  
On the other hand, turning our attention towards the posterior $P$, an important class of models is exponential families, and in particular  Gibbs distributions such as \VIZ{Gibbs}, also referred in the machine learning literature  as Boltzmann machines, \cite{Murphy:2012,MacKay:2003}. 
For such posteriors,  the minimization in \VIZ{VP1} yields an equivalent  variational free energy, 
\cite{MacKay:2003}. 
Indeed this observation was also made in the coarse-graining literature and was the starting point 
of coarse-grained parameterizations based on a similar principle to the variational inference 
in (\ref{VP1}), \cite{Shell2008, shell3,zabaras,NoidReview:2011}, and \cite{KP2013}.   
However, an additional complexity to (\ref{VP1}) arises in the coarse-graining case, 
where the coarse-graining map $\COP$ in (\ref{CG-map1})
enters in the optimization problem (\ref{ENT1}),
\begin{equation}\label{VP-CG1}
   \min_{\theta\in \Theta} \RELENT{\PUSHFWD{P}}{\bar Q^\theta} \COMMA
\end{equation}
where $\bar Q^\theta$ is the parametrized coarse-grained model and $P$ the microscopic Gibbs distribution.
In addition to the variational inference point of view that 
allows for optimal parameterization of a coarse-grained model
we  can also use \VIZ{VP-CG1} as means to compare 
different coarse-graining maps $\COP^{(i)}$. 
That is we can deploy \VIZ{VP-CG1} as an {\em Information Criterion} to assess 
and order  their relative effectiveness, i.e., $\COP^{(1)}$ is  a superior  coarse-graining 
map to $\COP^{(2)}$,  if and only if,
\begin{equation}\label{VP-CG2}
  \min_{\theta\in \Theta} \RELENT{\PUSHFWDA{P}}{\bar Q^\theta}\,<\,\min_{\theta\in \Theta} \RELENT{\PUSHFWDB{P}}{\bar Q^\theta} \, .
\end{equation}
The heuristic interpretation of \VIZ{VP-CG2} is that $\COP^{(1)}$ provides better 
information compression than $\COP^{(2)}$.
  Authors  in \cite{Foley} compare  different coarse-graining maps (resolutions)  based on the entropic component of the many body potential of mean force.

\subsection{Computational variational inference methods in coarse-graining.}

In  equilibrium systems where the coarse-graining of Gibbs distributions $P$  is considered, the coarse-grained family 
$Q^\theta$ in (\ref{VP-CG1}) is not necessarily mean field so there is no available analytic solution such as the ones discussed in \cite{MacKay:2003}. In this case the minimization (\ref{VP-CG1})
is carried out using numerical optimization methods, \cite{shell1,zabaras}.
Although one can consider steepest descent and Newton-Raphson type methods for this optimization problem, the most efficient methods are 
stochastic optimization methods in the spirit of the Robbins-Monro algorithm, for the latter see \cite{Bottou:2004} and references therein. 
In particular, in \cite{zabaras} the authors propose an algorithm that is essentially a stochastic optimization version of the Newton-Raphson algorithm. 
This method improves the Robbins-Monro algorithm by employing a {\em natural gradient} time-stepping, \cite{Amari:1998}. 
The natural gradient time-stepping  arises as part of the Newton-Raphson scheme for (\ref{VP-CG1}), 
since the Hessian of the relative entropy is exactly the Fisher Information Matrix.
A similar algorithm was also  proposed in machine learning for stochastic variational  inference,  \cite{Hoffman:2013}. 
Finally,
a Newton-Raphson method was introduced in \cite{KP2013} for the minimization of path-space relative entropy for dynamics and non-equilibrium systems discussed in Section 4.
There the role of the natural gradient is played by the {\em  path-space  Fisher Information Matrix}, which is a type of Fisher Information for dynamics introduced first  in the context of sensitivity analysis in  \cite{PK2013}.

\section{  Microscopic, coarse-grained  and reconstructed processes for molecular systems}\label{AtCG}
 
In this section we  describe a prototypical molecular system  at the   microscopic  scale and introduce its  coarse graining as a configuration transformation that lumps together degrees of freedom.  We define the 
underlying microscopic evolution in terms of  a general stochastic  differential equation and propose   
parametrized stochastic dynamics  as  a Markovian approximation of the coarse-grained process.   
We introduce a reconstruction  of the coarse process   defining a process that reintroduces the lost degrees 
of freedom and is approximating the microscopic evolution.  
The reconstructed process  serves  as  an  auxiliary process that connects the coarse-grained approximating 
dynamics with the microscopic dynamics on the same space. 
%
%
%
%
\subsection{ Microscopic  dynamics}\label{MicroDyn} We assume a prototypical system of $N$ (classical) 
molecules in a box of the fixed volume $V$ at the temperature $T$.
Let $q=(q_1,\dots,q_N) \in \R^{3N}$ describe the position vectors of the $N$ particles in 
the   microscopic description and $ p=(p_1,\dots,p_N) \in \R^{3N}$ the momentum vectors. 
We denote by $x = (q,p)\in \R^{6N}$ the joint vector of position and momentum.
We consider the evolution of the $N$ particles described by  a diffusion process   
$\PXT$, a continuous time Markov process  satisfying the stochastic differential equation  (SDE)
\begin{equation}\label{gSDE1}
  \begin{cases}
    dX_t = b(X_t)dt + \sigma(X_t) dB_t\COMMA  \;\; t\ge 0\COMMA  \\
    X_0\sim \mu_0\COMMA
  \end{cases}
\end{equation}
where $b(x)\in \R^{6N}$ and $\sigma(x) \in \R^{6N \times k},\ k\le 6N$ are the drift and diffusion 
coefficients and $B_t$  denotes the standard $k$-dimensional Brownian motion. 
 The notation   $X_0 \sim
\mu_0$ means that the random variable $X_0$ is distributed according to  the probability $\mu_0$.  
Throughout this paper we assume that the vector field $b(x)$ and 
the diffusion coefficient field $\sigma(x)  $ are such that the system
\VIZ{gSDE1} has a unique solution for all $t\geq 0$, \cite{oksendal2003}.  
The general form \VIZ{gSDE1} also accommodates the case of Hamiltonian equations of motions with
the Langevin thermostat used for describing molecular systems at equilibrium in a thermal bath. 
 
In general  
a Langevin process does not necessarily possess a stationary distribution, or  such distribution may not be explicitly known. 
 This can be the case  when  non-conservative external forces appear  or when  the  detailed balance condition fails. 
Note that for stationary dynamics
the detailed balance condition and time-reversibility  
are equivalent. 
If the force  $F(q) $ is conservative  and the fluctuation dissipation relation is satisfied,  
  the Langevin dynamics are time irreversible  
up to momentum reversal, \cite{Lelievre}, which 
guarantees that the Gibbs canonical measure is a stationary distribution
\begin{equation}\label{Gibbs}
 \mu(dq) =Z^{-1}\exp\{ -\beta \IP(q)\}dq \COMMA
\end{equation}
where $\IP(q) $ is the potential energy such that $F(q) =- \nabla\IP(q)$,    $Z= \int_{\R^{3N}}e^{-\beta \IP(q)} dq$ is  the partition function and $\beta=\frac{1}{k_B T}$,    $k_B$ the Boltzmann constant.
On the other hand for open systems, where the force  $F(q)$ has a non  conservative part,   \cite{Gallavotti1995}, i.e. 
$$F(q) \neq - \nabla\IP(q)\PERIOD$$ 
a  {\it non equilibrium steady state} (NESS)  exists  for which the condition of detailed balance fails.  
We note here that in machine learning and neural networks  distributions such as \VIZ{Gibbs}  are called Boltzmann machines and related reversible dynamics are considered for instance in stochastic Hopfield Models, \cite{MacKay:2003}.

%
%
%
%
\subsection{Coarse-grained   and reconstructed dynamics}\label{CGdyn}
Coarse-graining is considered as the application of a {\it linear} mapping (CG mapping) 
(see also (\ref{CG-map1}) for a more general definition)    
\begin{eqnarray*}
\COP:\R^{n} \to \R^{m}, \qquad 
 x  \mapsto \COP x   \in \R^{m} 
\end{eqnarray*}    
on the microscopic state space. For notational simplicity we set $n=6N$ and $m=6M$. 
The mapping  determines the $M(<N)$ CG particles with the state $\bx =\COP x$ as a function   
of the microscopic state $x$. 
Examples of CG maps commonly used  for molecular systems include the mapping to the centers of mass 
of groups of atoms, the end-to-end vector of molecular chains, projections to a collection of atoms, 
see also Section~\ref{Center}.
We call `particles'  
and `CG particles' the elements of the  microscopic and coarse  configuration space respectively.
 
The {\it proposed coarse space dynamics} are described by a Markov process $\{\bX_t\}_{t\ge 0} $ in $\R^m$    approximating the process $ \{\COP X_t\}_{t\ge 0}$ which is, in principle, non-Markovian.  
The Markov process  $\{\bX_t \}_{t\ge 0}$ is  given as the  solution of  the   {\it parametrized} stochastic differential equations 
\begin{equation}\label{CGSDE}
\begin{cases}
 d\bX_t = \bbr(\bX_t;\theta)dt + \BARIT\sigma(\bX_t;\theta) d\bB_t,  \ t > 0, \\
  \bX_0\sim \BARIT\mu_0\COMMA
 \end{cases}
\end{equation}
where   the drift $\ \bbr(x;\theta)\in \R^m $  and  diffusion $ \BARIT\sigma(x;\theta) \in \R^{m\times l}$, $\    l\le m$,  coefficients are parametrized with   $\ \theta \in \Theta$.
 $ \bB_t$ is an $l$-dimensional standard Brownian motion.
 As we have already indicated the goal of our study is to find the most effective among the proposed CG models 
    such that   $\{\bX_t \}_{t\ge 0}$  ``best approximates" the  process $ \{\COP X_t\}_{t\ge 0}$,
  that is  to  find  optimal   $\bbr(x;\theta) $ and  $\BARIT\sigma(x;\theta)$ in a parametric or non-parametric form, 
  which is the subject of Section~\ref{RERminS}. 
  
We define a {\it reconstructed } process 
of  a coarse process $\{\bX_t\}_{t\ge 0}$ in $\R^{m}$ onto the microscopic space $\R^{n}$ to be any stochastic process $\tPXT $ which satisfies
\begin{equation}\label{backmaprel}
 \COP \tX_t = \bX_t\COMMA \quad  \ t\ge 0 \text{ in distribution}\PERIOD
\end{equation}
A reconstructed process $\tPXT $  is  obviously not unique,  a trivial example is when $\COP$ is not a one-to-one transformation  though the non-uniqueness is not of a concern for our methodology.  
The path-space measures in $[0,T]$ of the process we consider are denoted by
\begin{eqnarray}\label{PathMeas}
&&\PATHPT \text{ for   the microscopic process } \PXT\COMMA \nonumber \\
&&\CGPATHPT \text{ for   the coarse  process } \CGPXT  \text{ and}\\
&&\tPATHPT = \PULLBCK{\bar Q^\theta}  \text{ for  the reconstructed process } \tPXT\COMMA \nonumber 
\end{eqnarray}
where the notation $  \PULLBCK{\bar Q^\theta} $ is described in Section~\ref{EqRE} as the reconstruction of the coarse path space measure $\CGPATHPT$ onto the microscopic space.   For the purpose of the present work  
we assume that the reconstructed process is the solution of the system
\begin{equation}\label{rSDE1}
 \begin{cases}
   d\tX_t = \tb(\tX_t;\theta)dt +  \tilde \sigma(\tX_t;\theta) dB_t\COMMA\;\;  \ t\ge 0\COMMA\\
   \tX_0\sim \nu_0 \PERIOD
 \end{cases}
\end{equation}
The coefficients   $\tb(x;\theta)$ and $ \tilde \sigma(x;\theta) $ must be such that 
the relation~\VIZ{backmaprel} is satisfied.  Note that from the definition of the reconstructed process we have that for every  observable of the form 
$$
f(x)=g(\COP x)\COMMA
$$
i.e., a coarse observable, the expectations with respect to the probability of the reconstructed and the coarse process are  identical 
\begin{equation}\label{CGobserv}
\E^x[f(\tX_\tau)] = \E^{\COP x}[g(\bX_\tau)] \COMMA
\end{equation}
 when   $\tX_0= x $ and $\bX_0=\COP x $ and   for  any stopping time $\tau>0$. We denote $\E^x$ the expectation with respect to the probability of $\{\tX_t\}_{t\ge 0}$, 
$R(\tX_{t_1}\in F_1\, \dots, \tX_{t_k}\in F_k) = P(\tX^x_{t_1}\in F_1\, \dots, \tX^x_{t_k}\in F_k)$,   for any $F_i, i=1,\dots,k $  subsets of $\R^n$,  where $\tX^x_{t}$ denotes that $\tX_t$ starts at $X_0=x$.  $\E^{\COP x}$ denotes the expectation with respect to the probability of   $\{\bX_t\}_{t \ge 0}$ starting at $\bX_0=\COP x$.

With the following theorem we give   sufficient  conditions  that the 
$\tb(x;\theta)$ and $ \tilde \sigma(x;\theta) $  must satisfy in order the relation~\VIZ{backmaprel} to hold.   
A detailed study of reconstructed processes is presented in the work \cite{D57}, where explicit examples for stochastic lattice systems are presented. The proof of the theorem is given in \ref{appendixProof}, 
and follows from the martingale uniqueness theorem, \cite{karatzas1991,oksendal2003}.

\begin{theorem}\label{PropRecon} {(\bf Reconstruction)}\\
 Let the processes  $\{\bX_t\}_{t\ge 0}$ in $\R^{m}$ and $ \tPXT $ in $\R^n$ be   solutions of 
 \VIZ{CGSDE} and \VIZ{rSDE1} respectively, and $\COP:\R^n\to \R^m$  a linear mapping. Assume that 
$\tb(x;\theta)$ and $ \tilde \sigma(x;\theta)$ are such that the existence and uniqueness of solutions to \VIZ{rSDE1} is guaranteed. 
If 
 \begin{eqnarray}
&&\COP  \tb(x;\theta) = \bbr(\COP x;\theta) \label{backdrift}\COMMA\\ 
&&\COP  \tilde \Sigma(x;\theta)  \COP^{tr} = \BARIT\Sigma(\COP x;\theta)  \ \text{ for all } x \in \R^n \COMMA \label{backdiff}
\end{eqnarray} 
  where $\tilde \Sigma(x;\theta) =  \tilde \sigma(x;\theta) \tilde \sigma^{tr}(x;\theta) $ and $\BARIT \Sigma(\bx;\theta) =   \BARIT\sigma(\bx;\theta) \BARIT\sigma^{tr}(\bx;\theta)$   and $\cdot^{tr}$ denotes matrix  transpose,
then 
  \begin{equation*}  
  \COP \tX_t = \bX_t , \ t\ge 0  \text{ in distribution}  \PERIOD 
   \end{equation*}
   In particular, the relation \VIZ{CGobserv} holds for all coarse observables $f(x) =g(\COP x)$.
\end{theorem}

 \begin{remark}\label{RemarkSec3}
 {\rm 
 (a) For a non-linear   CG mapping   $\COP$,  such as transformations to reaction coordinates \cite{KHKP1,Lelievre2010}, relations~\VIZ{backdrift} and \VIZ{backdiff}   do not hold.  The proof of Theorem~\ref{PropRecon} demonstrates the direction how these relations  can be generalized to non-linear mappings.   \\
 (b) In the rest of our work we assume that the diffusion term $ \BARIT\sigma(\bx;\theta) $ is not parametrized i.e., $ \BARIT\sigma(\bx;\theta)= \BARIT\sigma(\bx )  $  
    and is such that 
    \begin{equation*}\label{CGdiffus}
 \BARIT\Sigma(\COP x) =   \COP  \Sigma(x)  \COP^{tr}      \COMMA
   \end{equation*}
   where $  \Sigma(x) = \sigma(x)\sigma^{tr}(x)$ and $\BARIT \Sigma(\bx;\theta) =   \BARIT\sigma(\bx) \BARIT\sigma^{tr}(\bx)$.
This requirement ensures that the diffusion coefficient of the reconstructed process $\tPXT$, coincides with the one  of the microscopic process $\PXT$,  
  \begin{equation}\label{Eqdiff}
   \tilde\sigma(x;\theta ) =  \sigma(x)\COMMA 
   \end{equation}
a condition that we need for the development of the theoretical tools in Section~\ref{RERmin-sec2}.
We point out that the assumption \VIZ{Eqdiff} is not necessary in the development of the variational inference problem for data driven systems; this will be further discussed in Section~\ref{data_driven}.
}
\end{remark}

\section{Variational inference for coarse-grained dynamics.} \label{Sec:VarInf}
In analogy to the previous discussion in Section~\ref{EqRE} for the coarse-graining of Gibbs distributions,
our approach to coarse-graining of dynamics can also be viewed as a variational inference method, however,
this time set in the path space. 
We present an extension of the information-theoretic approach to systems with {\it non-equilibrium steady states} as well as dynamics in  {\em finite times}. The presented method also allows for approximation of {\it dynamical} observables, i.e., quantities that are averaged over the path distribution instead of over a distribution at a terminal time.
  
The relative entropy between two path measures $\PATHS$ and $\PATHSAPP$ (see \cite{KP2013} for a specific examples)
for the processes on the interval $[0,T]$ is 
\begin{equation}\label{path-relentropy}
\RELENT{\PATHS}{\PATHSAPP} = \int \log\frac{d\PATHS}{d\PATHSAPP}\, d\PATHS \equiv \EXPECT_{\PATHS}\left[\log\frac{d\PATHS}{d\PATHSAPP}\right]\COMMA
\end{equation}
where $\tfrac{d\PATHS}{d\PATHSAPP}$ is the  likelihood ratio  (or the Radon-Nikod\'ym derivative) of $\PATHS$ with respect to $\PATHSAPP$ and 
   $  \E_{P_{[0,T]} }[f] $ denotes averaging  over the probability of paths  $X_t(\omega)$, ${t\in [0,T]}$, $ P_{[0,T]}$. 
Naturally, we have to assume that the measures $\PATHS$, $\PATHSAPP$ are absolutely continuous,  that is if for an event $A$ holds $\PATHSAPP(A)=0$ then $\PATHS(A) =0 $, and $\log\frac{d\PATHS}{d\PATHSAPP}$
is $\PATHS$-integrable. 

In the setting of coarse-graining or model-reduction the measure $\PATHS$ is associated with the exact process (mapped to the coarse space) and 
$\PATHSAPP$ is associated with the approximating coarse-grained process. 
In general the relative entropy \VIZ{path-relentropy} in this dynamic setting is not a suitable object for analyzing steady states
or long-time behavior, however, 
in practically relevant cases of {\it stationary} Markov processes we can work with the {\em relative entropy rate} (RER) 
\begin{equation}\label{rate-relentropy}
\ENTRATE{P}{Q} = \lim_{T\to \infty} \frac{1}{T} \RELENT{\PATHS}{\PATHSAPP} \COMMA
\end{equation}
where $P$ and $Q$ denote the corresponding stationary processes. 
The definitions of (\ref{path-relentropy}) and (\ref{rate-relentropy})
do not depend on knowing the distributions or an underlying Gibbs distribution 
making them {\em suitable for non-equilibrium problems} where even the steady states are not known analytically. 
We also refer to \cite{PK2013,Tsourtis:2015,DKPP}  where this feature is employed to develop sensitivity 
analysis methods for non-equilibrium systems.

In order to select the best approximation of the coarse-grained model, i.e., 
of the exact coarse-grained measure $\PATHSAPP$, we define
a parametrized family of measures $\PATHSAPP^\theta$
depending on parameters $\theta\in \Theta$. The best approximation is fitted using entropy
based criteria in order to find the best Markovian approximation of the coarse-grained process. 
We consider the optimization principle
 \begin{equation}\label{VP-path1}
 \min_{\theta\in \Theta} \RELENT{P_{[0, T]}}{\PULLBCK{\bar Q^\theta_{[0, T]}}}\, ,
 \end{equation}
where 
$P_{[0,T]}$ is the path distribution of the original microscopic process and 
$\PULLBCK{\bar Q^\theta_{[0, T]}}$ 
is the parametrized path-space coarse-grained 
distribution
back-mapped to the microscopic space.
Furthermore, 
for coarse-graining of stationary dynamics 
we consider the variational inference optimization problem based on the relative entropy rate (RER)  instead of the full relative entropy
\begin{equation}\label{VP-path2}
  \min_{\theta\in \Theta} \ENTRATE{P}{\PULLBCK{\bar Q^\theta}}\PERIOD
\end{equation}
We note that  for dynamics 
we have essentially another Gibbs structure such as (\ref{Gibbs}),
however, this time in the space-time. Precisely this structure is used to obtain the relative entropy rate calculations in Theorem~\ref{thmRel}.

%
%
%
%
\subsection{Path-space information methods for diffusion processes}\label{RERmin-sec2}
The purpose of this section is  to   compare   continuous time Markov processes at the same state space given as  solutions of stochastic differential equations.
We provide explicit  representations of the  relative entropy (RE) and relative entropy rate (RER) for 
the path-space measures  of the  processes   $\PXT$ and  $\tPXT$   in $\R^n$,  
solutions of  the SDEs \VIZ{gSDE1} and  \VIZ{rSDE1} respectively, in terms of their drift and diffusion coefficients, for   (a) the finite time and (b) the stationary  (steady-state)   regime 
for molecular systems under non-equilibrium conditions.  
 
The key mathematical tool  that allows us to express the RE $\RELENT{\PATHPT}{\tPATHPT}$,  and the RER  $\ENTRATE{\PATHP}{\tPATHP}$, defined in \VIZ{path-relentropy} and \VIZ{rate-relentropy} respectively, in terms of the drift and diffusion coefficients appearing  in \VIZ{gSDE1} and \VIZ{rSDE1}, is the Girsanov theorem, see \ref{appendixProof2},  \cite{oksendal2003, karatzas1991}. 
Suppose that there exists a process $\{u(X_s;\theta)\}_{s\ge 0} $   in $\R^{k}$ such that
\begin{equation}\label{Ur}
\sigma(X_s)u(X_s;\theta) = b(X_s) -  \tb(X_s;\theta), \text{ and } \E\left[ \EXP{\frac12\int_0^T |u(X_s;\theta)|^2\,ds}\right] <\infty \COMMA
\end{equation}
where  $|u(x;\theta)|^2 =\sum_{i=1}^k  u^2_i(x;\theta)$.
Recall that a  process $\PXT$ is   stationary when its joint distribution 
does not change with time. 

\begin{theorem}\label{thmRel}
Let $\PXT$, $\tPXT$   be    Markov processes,  solutions of \VIZ{gSDE1} and \VIZ{rSDE1}   with path space distributions  $\PATHPT,\tPATHPT$, respectively, and  $ X_0 \sim \mu_0$, $\tilde X_0 \sim \nu_0$ for which \VIZ{Eqdiff} holds. Suppose that there exists  a process $\{u(X_t;\theta)\}_{t\ge 0}$, as  defined in \VIZ{Ur}. Then
\begin{enumerate}[a)]
\item the relative entropy between $\PATHPT$ and $\tPATHPT  $, for any $[0,T]$,  is  
\begin{equation}\label{RE}
\RELENT{\PATHPT}{\tPATHPT} =\Hh^T(\PATHPT\SEP\tPATHPT)  + \RELENT{\mu_0}{\nu_0}\COMMA
\end{equation}
where  
\begin{equation}\label{REFT}
 \Hh^T(\PATHPT\SEP\tPATHPT) =  \EXPECTWRT{\PATHPT}{\frac12 \int_0^T |u(X_s;\theta)|^2\,ds}  \PERIOD
\end{equation}
\item If furthermore $\PXT$ and  $\tPXT$ are  stationary Markov processes, with $\mu(dx)$   the invariant measure for $\PXT$ then
\begin{equation}\label{REst}
\RELENT{\PATHPT}{\tPATHPT} = T\ENTRATE{\PATHP}{\tPATHP}  + \RELENT{\mu}{\nu_0}\COMMA
\end{equation}
where $\ENTRATE{\PATHP}{\tPATHP}$ is the  the relative entropy rate  which is given by  
\begin{equation}\label{RERf0}
 \ENTRATE{\PATHP}{\tPATHP} =  \EXPECTWRT{\mu}{\frac12 |u(X;\theta)|^2} \PERIOD
\end{equation}
  \end{enumerate}
\end{theorem}

The proof of the theorem is given in \ref{appendixProof3}.
Theorem~\ref{thmRel}~a) gives a form of the relative entropy {\em for any finite time interval $[0,T]$} while Theorem~\ref{thmRel}~b)  addresses the  long time regimes providing a reduced form of \VIZ{RE}.
 Note that the relation \VIZ{RE} holds for any process $\PXT$, not necessarily stationary. 
 Another important result that Theorem~\ref{thmRel}  states is that relation \VIZ{RERf0} holds for any  initial condition $\tilde X_0 \sim \nu_0$, where $\nu_0$ is not necessarily an invariant measure.    The  two  properties that ensure  the linear in time  increase of the RE and the existence of the RER are (a) 
 the Markov property and (b) the stationarity. 
Thus the ratio of the path probability densities on $[0,T]$ will be constant independent of 
the  time interval  length $T$. 
Markovianity though is not the most general condition, for example the semi-Markov property is also sufficient, \cite{PK2013}.

As we  see from relation~\VIZ{RE}  the 
relative entropy $\RELENT{\PATHPT}{\tPATHPT}$ is a quantity that increases in time, thus calculation of the RE  over long time intervals may become unfeasible. 
Though for finite time intervals calculations may be tractable. 
This fact is depicted in the following corollary where  we provide an explicit formula of the
 finite time component \VIZ{REFT} along with one of RER \VIZ{RERf0} in terms of the drift and diffusion 
 coefficients.
 
\begin{corollary} \label{thm3}
\begin{enumerate}[(a)]
\item 
\label{thm3a} 
{\bf (RER representation for stationary processes)}\\
Let  $ \sigma(x) \in \R^{n\times k}, \ x\in \R^n$  appearing in \VIZ{gSDE1} with $\text{rank}(\sigma(x)) = r = k\le n$.
Suppose that there exists $u(x;\theta)\in\R^k $ satisfying  \VIZ{Ur}. Then
\begin{equation}\label{RERf1}
          \ENTRATE{\PATHP}{\tPATHP}  
           =  \EXPECTWRT{\mu}{\frac12 \|  b(X) -  \tb(X;\theta)\|^2_{\Xi}}   \COMMA
\end{equation}
  where $\Xi(x) = \left[\sigma^{tr}(x)\sigma(x)\right]^{-1}\sigma^{tr}(x)$ 
 and  $\| \cdot \|_\Xi $ denotes  the norm
\begin{equation}\label{XiNorm}
 \| z\|_\Xi 
 =  z^{tr}\Xi^{tr} \Xi z\, ,\  z \in \R^n\PERIOD
\end{equation}
  \item \label{FTRE} {\bf(RE representation for finite time)}
\begin{equation*} 
\RELENT{\PATHPT}{\tPATHPT} =\Hh^T(\PATHPT\SEP\tPATHPT)  + \RELENT{\mu_0}{\nu_0}\COMMA
\end{equation*} 
where
\begin{equation}\label{RET}
 \Hh^T(\PATHPT\SEP\tPATHPT) =  
 \EXPECTWRT{\PATHPT}{\frac12 \int_0^T  \|  b(X_s) -  \tb(X_s;\theta)\|^2_{\Xi}   \,ds}  \PERIOD
\end{equation}
\end{enumerate}
\end{corollary}

The result of Corollary~\ref{thm3} can be generalized to any $\sigma(x) $ 
with $\text{rank}(\sigma(x)) = r< k\le n$ if we use in place of $\Xi(x)$ 
a (Moore-Penrose) generalized inverse  
of $\sigma(x)$, \cite{RomanALA, ben2003generalized}.
For completeness we provide the proof of Corollary~\ref{thm3}(a)  
in \ref{appendixProof4} and note that the proof of (b) is similar.

\begin{remark}
{\rm Note that the representations \VIZ{RERf1} and \VIZ{RET} are valid when $\PXT$ or   $\tPXT $ are Ito processes.
The first time we need that $\PXT$ and $\tPXT $ are  both (Ito) diffusions is when we want to use the stationarity of the processes to simplify the 
RE to RER  such  that we have $u(x;\theta)$ independent of time $t$.
If $b(x)$ is substituted by $b(x,t) $ (in which case the microscopic process is non-Markovian),  
\cite{oksendal2003},  
  then  we would have 
$$\Hh^T(\PATHPT\SEP\tPATHPT)= \EXPECTWRT{\PATHPT}{\frac12 \int_0^T  \|  b(X_s,s) -  \tb(X_s;\theta)\|^2_{\Xi}   \,ds}  \PERIOD$$
Even if $\PXT$ is stationary  it is not obvious that $  b(X_s,s) -  \tb(X_s;\theta)$ is stationary, and the reduction to the RER needs to be checked casewise.
 }
\end{remark}

\section{Variational inference in coarse graining and path-space force
matching}\label{RERminS}

The goal of this section is to obtain optimal parameters $\theta^*$ of a coarse-grained
dynamics model $\CGPXT$, eq. \VIZ{CGSDE}, approximating the microscopic dynamics $\PXT$, eq. \VIZ{gSDE1}, based on the path-space variational inference described by  \VIZ{VP-path1} and \VIZ{VP-path2}.
We prove that minimization of the RE reduces to  time dependent  weighted least squares type problems 
with weights that depend on the CG mapping and the diffusion coefficient. 
For stationary regimes the time dependence  is altered with the minimization of the RER. 
The weighted least squares formulation provides a natural generalization of the force matching 
methods developed for systems at equilibrium, \cite{IzVoth2005}, to  non-equilibrium systems. 
Moreover, the  relation of RE  and force matching type optimization methods is revealed, 
in analogy to the  equilibrium case, studied in  \cite{Voth2008a, NoidReview:2011, KHKP1}.

%
%
%
%
\subsection{Properties of the CG mapping and the reconstructed process}\label{BackMaps}
We recall that the coarse-graining   map 
$\COP:\R^n\to \R^m$ defined in Section~\ref{AtCG} is a {\it linear transformation}. 
We denote with the same letter $\COP\in \R^{m\times n}$ the matrix representation of the linear map $\COP$.
Furthermore, we assume that the CG map has full rank, 
\begin{equation*}
\text{rank}(\COP) = m\PERIOD
\end{equation*}

We consider  the reconstructed process $\tPXT$ as  defined in Section~\ref{CGdyn}, Theorem~\ref{PropRecon}, and  Remark~\ref{RemarkSec3}, i.e., such that  
 \begin{equation*}\label{brec1}
 \COP \tb(x;\theta) = \bbr(\COP x;\theta),\ \text{and } \tilde{\sigma}(x;\theta) = \sigma(x) \text{ for all } x \in \R^n\PERIOD
 \end{equation*}
  Thus the coarse diffusion coefficient $\BARIT \sigma( \bx )$  is  independent of the parameters $\theta$,  and
 \begin{equation}\label{CGdiff1}
  \BARIT\Sigma(\COP x ) = 
\COP   \Sigma(x )  \COP^{tr}    \ \text{ for all } x \in \R^n \PERIOD
 \end{equation}
 The reconstructed drift term $ \tb(x;\theta)$ is a  vector field in $\R^n$  that can be written as 
 \begin{equation}\label{brecRep}
  \tb(x;\theta) =\GCOP\bbr(\COP x;\theta) + \left( I_n -\GCOP \COP \right) y^{\perp}(x), \ \text{ for all } x \in \R^n \COMMA
 \end{equation}
where   $I_n$ denotes the identity matrix in $\R^n$, 
$\GCOP$ is a right inverse of $ \COP$, an $ n\times m $ matrix such that 
 \begin{equation}\label{CGLeftInverse}
 \COP\GCOP=I_m\COMMA
 \end{equation}
and $y^{\perp}(x)$ is any arbitrary vector in $\R^n$ satisfying $\COP\left( I_n -\GCOP\COP \right)y^{\perp} = 0_m $.
In principle, $\GCOP$ is not unique, though for any such $\GCOP$ it holds $\BARIT b =\COP \tb$
which is  the main property that we need. We choose $\COP$ to have   the full rank 
so  we   have an explicit form  
$\GCOP$,  that is 
\begin{equation*}\label{SpCOPInv}
  \GCOP=\COP^{tr} \left( \COP\COP^{tr}\right)^{-1}\PERIOD
  \end{equation*} 
Most of the CG maps related to specific applications we consider are of a full rank, e.g.,
mapping to the centers of mass of   groups of atoms, or projecting to fewer state space coordinates, 
see Section~\ref{Center}.
As already mentioned, and   relation \VIZ{brecRep} verifies,  the reconstructed process is not unique. 
Therefore one can always  choose   the  term of  the reconstructed drift  $\left( I_n -\GCOP \COP \right) y^{\perp}(x)$ in \VIZ{brecRep}  independent of the parameter $\theta$.
In the rest of this work we assume that   $\left( I_n -\GCOP \COP \right) y^{\perp}(x)$   is independent of the parameter $\theta$.

%
%
%
%
\subsection{Optimal coarse-grained dynamics} 
Having set up the optimization problems \VIZ{VP-path1} and \VIZ{VP-path2}, in this section we look for optimal solutions $\theta^*(T)$  and $\theta^*$    respectively,
based on the first order optimality condition. That is  if $\theta^* $ is a  solution  of \VIZ{VP-path2}   then 
\begin{equation}\label{Optimality}
\nabla_\theta \ENTRATE{\PATHP}{\tPATHPOP}  = 0\PERIOD
\end{equation}
Thus  solutions of \VIZ{Optimality} reveal the  local optima of the RER.
Note that if the RER  is a strictly convex function of $\theta$ then there is  
a unique (global minimum) $\theta^*$. This property clearly  depends  on the choice of 
the parametrized model, i.e., through the definition of the parametrized drift 
$\bbr(\bx; \theta)$, see for an example Remark~\ref{NonParamRemark}.

\begin{theorem}\label{REminThm}
 Let   $\COP:\R^n\to \R^m$ be a linear mapping with $ \text{rank}(\COP) = m$.
Consider  the microscopic process $\PXT$  and the coarse process $ \CGPXT$ satisfying ~\VIZ{gSDE1} and ~\VIZ{CGSDE} respectively, 
 with the drift $b(x) $, $\bbr(\bx;\theta) $ and the
  diffusion   terms $\sigma(x) $, $ \BARIT\sigma(\bx)$, such that  \VIZ{CGdiff1} holds,   
  $ \text{rank}(\sigma) = k$ 
  and  $X_0\sim \mu_0$, $\BARIT X_0\sim \BARIT\mu_0$.
Let $\tPXT$ be a reconstructed process of  $ \CGPXT$ with the 
drift $ \tb(x;\theta) $ defined in \VIZ{brecRep} and the diffusion coefficient $\sigma(x) $.
 Then, for any   $\GCOP $ satisfying \VIZ{CGLeftInverse},
 \begin{enumerate}[a)]
\item 
 \begin{equation*}
\mathrm{argmin}_{\theta\in\Theta}\RELENT{\PATHPT}{\tPATHPT} = \mathrm{argmin}_{\theta\in\Theta}  
  \EXPECTWRT{\PATHPT } { \frac12 \int_0^T \| \COP b(X_s) - \bbr(\COP X_s;\theta) \|_{\GCOP\Xi}^2\,ds } \PERIOD
 \end{equation*}
where  
\begin{equation*}\label{CGXiNorm}
\| z \|_{\GCOP\Xi}^2 = z^{tr} \COP^{\sharp,tr}\Xi^{tr} \Xi\GCOP z,\ z\in \R^m \text{ and }  \ \Xi  = \left[\sigma^{tr}(x)\sigma(x)\right]^{-1}\sigma^{tr}(x) \PERIOD
\end{equation*}
\item
If  moreover  $\PXT$ is  stationary with the invariant measure $\mu$, then
 \begin{equation*}
\mathrm{argmin}_{\theta\in\Theta}\ENTRATE{\PATHP}{\tPATHP} =  
\mathrm{argmin}_{\theta\in\Theta}\EXPECTWRT{\mu} {\frac12  \| \COP b(X) - \bbr(\COP X;\theta) \|_{\GCOP\Xi}^2 }\COMMA  
  \end{equation*}
\end{enumerate}
\end{theorem}
We present the proof in \ref{appendixProof5}.

Theorem~\ref{REminThm}
 proves that the variational inference problems \VIZ{VP-path1} and \VIZ{VP-path2}  
are force matching type problem with the norm $\| z\|_{\GCOP\Xi}$ instead of the usual Euclidean norm. 
Moreover,  Theorem~\ref{REminThm}~{a)}  is a time dependent force matching, i.e., force matching over  paths, 
thus we may call the optimization problem the  {\it path-space force matching}. 
Note that  the calculations involved in the proof of Theorem~\ref{REminThm}  are independent of the fluctuation-dissipation 
and detailed balance, thus they  apply to non-equilibrium systems directly. 
  Theorem~\ref{REminThm}~b)   shows that relative entropy rate (RER) minimization  and FM are essentially identical. 
In particular formula  \VIZ{REst}  mathematically explain 
the difference between RE and FM at equilibrium, it 
reveals  what each method is doing: the one minimizes   $\Rr(\mu|\mu^\theta)$
 and the other the relative entropy rate $H(P|Q^{\theta})$.

 \begin{remark}\label{NonParamRemark}
 {\rm
 The existence of the {\em unique global minimum}  is  guaranteed if,  for example, 
 the force field $ \bbr(\bx;\theta)$ depends linearly on $\theta$.
Let $ \{\phi_k(\BARIT x)\}_{k=1}^{K}$ be a set of $m$-valued polynomials on $\R^m$ and approximate 
the force field $ \bbr(\bx;\theta)$ by
\begin{equation*}\label{nonparametric}
\bbr(\BARIT x; \theta) = \sum^K_k \theta_k \phi_k(\BARIT x)\PERIOD
\end{equation*} 
In this case the minimization of RER  is the least squares fit with respect to the stationary measure $\mu(dx)$. Due to the linear dependence   on $\theta$ the minimization problem has a unique solution due to strict convexity of  RER. Note that RE is also convex in this case and uniqueness of the corresponding optimization problem is guaranteed. 
The optimization problem reduces to solution of
the linear  system
\begin{eqnarray*}\label{PolBasis}
\BPhi \theta = \ba, \quad \text{ where }\;\; \BPhi_{ij} =    \EXPECT_\mu[\langle \phi_i,\phi_j\rangle_{\Pi^\sharp\Xi}]\COMMA\;\; \ba_i =   \EXPECT_\mu[\langle \phi_i,b\rangle_{\Pi^\sharp\Xi}]\COMMA\;\;
i,j =1, \dots K\PERIOD
\end{eqnarray*}
 Moreover, this set-up  
   is also used in estimation of non-parametric models, for example by
 specifying the basis set $\{\phi_i\}$  to be splines, or wavelets, \cite{Danani2015}.}
\end{remark}

\section{Effective Langevin dynamics and path-space force matching.}\label{Langevin-sec}
%
%
%
%

We present the application  of the path-space force matching  described in the previous section   
to the stochastic Langevin dynamics for molecular systems.
We propose  Langevin-type coarse dynamics  with a parametrized force (and friction) and derive the optimal parameter set for which  
the path space relative entropy is minimized, both for finite time  and stationary  dynamics.  
Furthermore, we present two   specific coarse graining maps:
  (a) The  transformation to the centers of mass of groups of atoms, and 
  (b) the  projection to a selected subset of atoms. 
 
The Langevin dynamics
are described by   the process $\{(q_t,p_t)\}_{t\ge 0} $, for an $N$-particle molecular system 
 with positions $ q=(q_1,\dots,q_N)  \in \R^{3N}$  and momenta $ p=(p_1,\dots,p_N) \in \R^{3N}$  which  
  satisfies 
 \begin{equation}\label{Lang}
 \begin{cases}
  dq_t= \mass^{-1} p_t dt \COMMA\\
 dp_t= F(q_t) dt -  \gamma \mass^{-1}p_t dt +\sigma  dB_t\COMMA
 \end{cases}
 \end{equation}
a Hamiltonian system coupled with a thermostat,  where  $ F(q)$ is the force field that is not necessarily a gradient, see also discussion in Section~\ref{MicroDyn}. 
 $ \mass=\DIAG(m_1I_3,\dots,m_N I_3)\in \R^{3N\times 3N}$ is the mass matrix, $\gamma\in \R^{3N\times 3N}$  is the friction and 
$\sigma\in\R^{3N\times 3N}$ the diffusion  coefficients respectively, and $B_t$ is the 
$3N-$dimensional Brownian motion. The diffusion and friction coefficients satisfy
the fluctuation-dissipation relation $ \sigma\sigma^{tr} = 2\beta^{-1} \gamma$.

Note that we have assumed  that  friction and diffusion coefficients are independent of the process, though our methodology is also applicable for $\gamma= \gamma(q)$  and $\sigma=\sigma(q)$, see Remark~\ref{Par_friction}.
The latter can be crucial   if  the microscopic and the CG molecular systems are not diffusive or when one wants to restrict the thermostat at the boundaries.
   
Let $M(<N)$ be the degrees of freedom of the coarse space, i.e., the number of CG particles,  and  define  the linear 
coarse-graining map $\COP_q: \R^{3N} \to \R^{3M}$   by
\begin{equation}\label{linearMap}
\COP_q q_j = \sum_{i=1}^N \zeta_{ji} q_i, \quad j=1,\dots,M, \text{ for any }   q\in \R^{3N} \COMMA
\end{equation}
for any set $\zeta_{ji}\in \R$,  $j=1,\dots, M, i=1,\dots, N $,  such that $\mathrm{rank}(\COP_q) = 3M$, 
  $\COP_q$ denotes the matrix representing the transformation \VIZ{linearMap}
\begin{equation*}
\COP_q = \left[ \zeta_{ji} I_3 \right]\COMMA
\end{equation*}
where $I_3 $ is the 3-dimensional identity matrix.
Furthermore, we denote 
 $$
 \COP:\R^{6N} \to \R^{6M}, \quad  \COP x = (\COP_q q, \COP_p p),\ x=(q,p)\COMMA
 $$
 where $\COP_p:\R^{3N} \to \R^{3M}$ denotes the momentum transformation.  
 Let the CG particles have mass matrix
$$ \cgmass  = \DIAG(\BARIT m_1I_3,\dots,\BARIT m_M I_3)\in \R^{3M\times 3M}\PERIOD$$
The momentum mapping $\COP_p$ is given by
\begin{equation}\label{massCGmap}
\COP_p =\cgmass\COP_q \mass^{-1}\COMMA
\end{equation}
such that   $ \COP_q dq_t = \COP_q\mass^{-1} p_t dt= \cgmass^{-1} \COP_p p_t dt $.
In the work \cite{Voth2008a} the mass matrix $\cgmass$ is defined such that a consistency condition 
in the momentum space is satisfied which,  with \VIZ{massCGmap}, defines  the mass  of the CG particles 
$\BARIT m_j, j=1,\dots, M $.  
The consistency condition states  that  the  momentum probability distribution of the coarse variables 
is the same on the coarse space and the microscopic space, 
that is  
$ e^{-\beta   \bp^{tr} \cgmass^{-1} \bp} \propto \int e^{-\beta  p^{tr} \mass^{-1} p} \delta(\COP_p p - \bp)\,dp$.  

The proposed  dynamics for the coarse variables $\bx =(\bq,\bp)\in \R^{6M}$  are given by the Langevin system
\begin{equation}\label{CGLang}
 \begin{cases}
  d\bq_t=   \cgmass^{-1}\bp_t dt \COMMA\\
 d\bp_t= \BARIT F(\bq_t;\theta) dt -  \BARIT \gamma \cgmass^{-1} \bp_t dt + \BARIT \sigma d\bB_t\COMMA
 \end{cases}
\end{equation}
where $\bB_t$ is a $3M$-dimensional Brownian motion. 
The diffusion coefficient $\BARIT \sigma$ is defined, according to relation \VIZ{CGdiff1},  such that
  $$
  \BARIT \sigma \BARIT \sigma^{tr} = \COP _p\sigma\sigma^{tr}  \COP_p^{tr}\PERIOD
  $$  
We examine two cases for the friction coefficient:
\begin{enumerate}[ (a)]
  \item    It matches the coarse graining transformation of the microscopic  friction forces   
     \begin{equation}\label{friction1}
     \BARIT \gamma \COP_q = \COP_p\gamma \PERIOD
     \end{equation}
  \item  The fluctuation-dissipation relation   is   satisfied for the coarse space dynamics. 
         Since the fluctuation dissipation relation is satisfied for the microscopic dynamics the friction coefficient in CG dynamics must be  
  \begin{equation}\label{friction2}
  \BARIT \gamma = \frac12 \beta  \BARIT \sigma \BARIT \sigma^{tr} =\COP _p \gamma \COP_p^{tr}\PERIOD 
  \end{equation}
   \end{enumerate}
Hence, the parametrization of the coarse dynamics is described here  only through the force 
$$\BARIT F(\bq;\theta), \quad \theta\in \Theta\COMMA$$
where  $\Theta$ is the parameter space.
The Langevin system    \VIZ{Lang} is   written in the form   of  the SDE  \VIZ{gSDE1} 
if we set $n=6N$, $x= (q ,p )$, 
$$ 
b(x)= (\mass^{-1}p , F(q)  - \gamma \mass^{-1} p )^{tr} \COMMA
$$
 and 
 $$\sigma (x) = \sigma_{0} =(0_{3N},\sigma)^{tr} \in \R^{6N\times 3N}\PERIOD$$
We also associate the coarse-grained dynamics with  the SDE  \VIZ{CGSDE}  where   $m=6M$, $\bx= (\bq ,\bp )$ and
 $$ \bbr(\bx;\theta)= (\cgmass^{-1}\bp , \BARIT F(\bq;\theta)   - \BARIT \gamma \cgmass^{-1}\bp )^{tr} \PERIOD$$
 Note that $\text{rank}(\sigma_{0}) = 3N $ and     
 $\text{rank}(\COP) = 6M $   since $\text{rank}(\COP_q) = 3M $.
  Moreover, we assume that given the CG map we can find a reconstructed process  
as described  in Section~\ref{CGdyn}.   
 
Applying the results of Sections~\ref{RERmin-sec2} and~\ref{RERminS} we provide the optimal parameter set $\theta^*$ 
for which the CG Langevin process $(\bq_t, \bp_t)_{t\ge 0}$  best approximates 
 $\COP(q_t, p_t)_{t\ge 0}=(\COP_q q_t, \COP_p p_t)_{t\ge 0}$, in the sense that  the path space RE (or RER) is minimized, 
 see Theorem~\ref{REminThm}. We state the result for the stationary regime and for the finite-time evolution.

\medskip
\noindent 
{\bf Stationary regime.}
Considering the description of the CG dynamics and assumptions on  $ \BARIT{\gamma}$, $\BARIT \sigma$ and $\COP$ mentioned above, 
conditions of Theorem~\ref{REminThm}~b) are satisfied  
thus we have that, for stationary microscopic Langevin dynamics
 \begin{equation*}
 \theta^* = \text{argmin}_\theta \EXPECTWRT{\mu} {\frac12  \| \COP  b(X) - \bbr(\COP  X;\theta) \|_{2,\GCOP\Xi}^2 }
 \end{equation*}
 where $\mu(dq,dp)$ is the stationary distribution for $\{(q_t,p_t)\}_{t\ge} 0$ and 
 \begin{equation*}\label{Xi2Norm}
  \| z \|_{2,\GCOP\Xi}^2 = z^{tr} \COP^{\sharp,tr}\Xi^{(2),tr}\Xi^{(2)}\COP^{\sharp}  z, \ z\in \R^{6M} \COMMA
  \end{equation*}
 with 
$$
 \Xi^{(2)} 
 = [\sigma^{tr} \sigma]^{-1} \sigma_{0}^{tr}\COMMA
$$
  and $\COP^{\sharp} \in \R^{6N\times 6M}$ is a right inverse of $\COP$, that has the form
  $\COP^{ \sharp}= \begin{bmatrix}
\GCOP_q &0  \\
0 &  \GCOP_p
\end{bmatrix} $ 
for any    $\GCOP_q\in \R^{3N\times 3M} $ such that 
$  \COP_q\GCOP_q =I_{3M}$ and  $ \GCOP_p  = \mass \GCOP_q \cgmass^{-1}  $.

 Therefore, 
\begin{equation*} 
 \|\COP b(x) - \bbr(\COP x;\theta)\|_{2, \COP^\sharp\Xi}^2 =  
  \| \COP_p F(q) - 
 \BARIT F(\COP_q q;\theta) - ( \COP_p \gamma \mass^{-1} p   -  \BARIT \gamma \cgmass^{-1} \COP_p   p )\|_{\GCOP_p \Xi}^2 \COMMA
 \end{equation*}
 where    
\begin{equation*} 
 \| u \|_{\GCOP_p \Xi}^2 = u^{tr}\COP_p^{\sharp,tr}\Xi^{tr}\Xi\GCOP_p u\, ,  \quad u\in \R^{3M} 
 \end{equation*}
and
 \begin{equation*}\label{XiPi}
 \Xi =[\sigma^{tr} \sigma]^{-1} \sigma^{tr} \PERIOD
\end{equation*}
If the friction coefficient is given as in  the case (a), eq. \VIZ{friction1}, 
the optimal parameter set is  
 \begin{equation}\label{Opt1}
 \theta^* = \text{argmin}_\theta \EXPECTWRT{\mu_p} {\frac12  \| \COP_p F(q) - \BARIT F(\COP_q q;\theta) \|_{\GCOP_p\Xi}^2 }     \PERIOD
 \end{equation}
 This is exactly the force matching problem related to the Hamiltonian dynamics.
On the other hand, if  the friction coefficient is given by (b), eq.  \VIZ{friction2}, 
  \begin{equation}\label{Opt2}
 \theta^* = \text{argmin}_\theta \EXPECTWRT{\mu_p} {\frac12  \| \COP_p F(q) - \BARIT F(\COP_q q;\theta)- ( \COP_p \gamma   -   \COP _p \gamma \COP_p^{tr} \COP_q ) \mass^{-1} p  \|_{\GCOP_p\Xi}^2 }  \PERIOD
 \end{equation}
 \begin{remark} 
 {\rm
The appearance of the term  $( \COP_p \gamma   -   \COP _p \gamma \COP_p^{tr} \COP_q ) \mass^{-1} p $ in the optimization is the result of the difference of the  friction forces that contribute to the CG particle, since we consider a fixed diffusion term $ \BARIT \sigma$, which is related to the  chosen stochastic CG dynamics \VIZ{CGLang}.  
}
 \end{remark}
 \begin{remark} 
 {\rm
Note that if    $ F(q) = - \nabla U(q)$,  then 
    \begin{equation*}\label{qpGibbs}
     \mu_p(dq,dp) = Z^{-1} \exp\{-\beta H(q,p)\}\COMMA
     \end{equation*}
       is the Gibbs distribution, where     $H(q,p) = \frac12 p^{tr} \mass^{-1} p +U(q)$.
If   $ F(q) \neq - \nabla U(q)$, i.e., driven Langevin dynamics, the form of  
$ \mu_p(dq,dp)$ is not analytically known, though in numerical implementations it is not needed as  
sets of samples $ (q^{(i)},p^{(i)})$ used to estimate the average are found as solutions of the Langevin system
at the stationary regime.
}
 \end{remark}

\medskip 
 \noindent
{\bf Finite time regime.} 
For parametrization on a finite time interval $[0,T]$, 
where the system has not reached stationarity   the approach to find (time dependent)
optimal parameters  $\theta(T)$ is given by application of  Theorem~\ref{REminThm}(a). 
 Using the same steps as for the stationary case,  we have that the optimal parameter  set for which the  process  $\{(\bq_t,\bp_t)\}_{t\ge} 0$ best approximates  $\{(q_t,p_t)\}_{t\ge} 0$ at the time interval $[0,T]$ is given  by, if  $  \BARIT \gamma  \COP_q=\COP_p\gamma$,
\begin{equation}\label{Opt3}
 \theta^*(T) = \text{argmin}_\theta \EXPECTWRT{\PATHPT } { \frac12 \int_0^T  \| \COP_p F(q_s) - \BARIT F(\COP_q q_s;\theta) \|_{\GCOP_p\Xi}^2 ds}  \COMMA
 \end{equation}
  and if   $\BARIT \gamma  = \tfrac{1}{2}\beta  \BARIT \sigma \BARIT \sigma^{tr} $
  \begin{equation}\label{Opt4}
 \theta^*(T) = \mathrm{argmin}_\theta \EXPECTWRT{\PATHPT } { \frac12 \int_0^T  
 \| \COP_p F(q) - \BARIT F(\COP_q q_s;\theta)- ( \COP_p \gamma   -  \COP _p \gamma \COP_p^{tr} \COP_q ) \mass^{-1} p_s  \|_{\GCOP_p\Xi}^2 ds} \PERIOD
 \end{equation}
 
 \begin{remark}\label{Equilibrium}
 {\rm {\bf (Connection to equilibrium force matching) }
  When the system is at  equilibrium   the optimization problem is defined by \VIZ{Opt1} and
   \VIZ{Opt2}, where the expectation is taken  with respect to the Gibbs measure $ \mu_p(dq,dp)$.
   As the averaged quantity in   \VIZ{Opt1} does not depend on $p$ the   $ \mu_p(dq,dp)$ expectation is equal  to the  expectation with respect to  $ \mu (dq)= Z^{-1} \exp\{-\beta U(q)\}\PERIOD $ 
  Therefore the optimization problem is 
   \begin{equation*} 
 \theta^* = \mathrm{argmin}_\theta \EXPECTWRT{\mu } {\frac12  \| \COP_p F(q) - \BARIT F(\COP_q q;\theta) \|_{\GCOP_p\Xi}^2 }     \COMMA
 \end{equation*}
  where $\COP_p=\cgmass \COP \mass^{-1}$ and $\| u\|^2 = u^{tr}\COP^{\sharp,tr}_p\Sigma^{-1}\GCOP u $, 
 with $\Sigma=\sigma\sigma^{tr}$ the covariance matrix. 
 Moreover, if the diffusion coefficient $\sigma $ is constant, identical for all particles, then $\Sigma^{-1}=\sigma^{-2}I_{3N}$ and if  $\COP_p^{\sharp,tr} \GCOP_p =I_{3M}  $, then
  \begin{equation*} 
 \theta^* = \text{argmin}_\theta \EXPECTWRT{\mu } {\frac{1}{2\sigma^2}  \| \COP_pF(q) - \BARIT F(\COP_q q;\theta) \|^2 }     \COMMA
 \end{equation*}
 where $\|\cdot\|$ denotes the  Euclidean norm in $\R^{3M}$. 
Therefore the path-space force matching is equivalent to the force matching 
for systems at equilibrium. 
  Note that $\COP_q^{\sharp,tr} \GCOP_q =I_{3M}$ holds, for example, if  $\COP_q$ is the mapping 
 to the group's center of mass  or an  orthogonal map, see the examples in Section~\ref{Center}.
 }
 \end{remark}
 \begin{remark}
 {\rm 
We should also state here that in all above applications a typical diffusive process is assumed for the dynamics of
the molecular systems. 
However, for realistic molecular systems this might not be a good approximation due to sub-diffusive behavior 
and the importance of long-time dynamical behavior or memory in a generalized Langevin equation (GLE)  framework. 
In such a case, either a scaling of the CG dynamics in a post-processing stage, or an approximation of 
the memory terms appeared in GLE through a parameterization of the diffusion term in the Langevin equation 
are required \cite{Ottinger,Harmandaris2009a,Harmandaris2009b,Fritz2011a,Buscalioni2010,Guenza2011,Briels}. 
This is the topic of a future work. 
}
\end{remark}
  
\begin{remark}\label{Par_friction}
{\rm {\bf (Parametrized friction)}
Note that if we consider that $\BARIT \gamma = \BARIT \gamma(q;\vartheta)$ then the same methodology can be applied with 
$\bbr(\bx;\theta,\vartheta) =  (\cgmass^{-1} \bp , \BARIT F(\bq;\theta)   - \BARIT \gamma(q;\vartheta) \cgmass^{-1}\bp )^{tr}\PERIOD$
For this generalization one should be careful to ensure the existence of  
CG processes with such friction and  of its stationary measure. For example, consider whether fluctuation-dissipation condition holds.  
We also refer to related works \cite{Briels2001, VothDyn2006, Buscalioni2010, Karniadakis2015}. 
}
 \end{remark} 
\subsection{Numerical implementation and optimal parameters}
 In \ref{numerics-sec5} we study the RER minimization problem
induced  by discrete time numerical schemes of the  Langevin dynamics  and overdamped Langevin dynamics. 
We consider the Langevin   \VIZ{Lang} and the coarse space   Langevin \VIZ{CGLang} dynamics,   apply the Brooks-Brunger-Karplus   (BBK), \cite{BBK1984},     discretization  scheme  in both  systems    and  study the RER  
  defined for the discrete  Markov chains,  see \VIZ{Path-L7}. 
 
We also study  
  the discrete time  RER minimization problem 
for the  overdamped Langevin dynamics $\PXT$ in $\R^n$,  
\begin{equation}\label{ovdLang}
dX_t = -\frac12 \Sigma(X_t)\nabla U (X_t) dt +  \frac12 \nabla  \Sigma(X_t) dt +\sigma(X_t) dB_t\COMMA
\end{equation}
  with a  $\theta$-parametrized coarse-grained approximation considered in the spirit of Section~\ref{CGdyn},   
  presented in detail in \ref{numerics-sec5}, using the Euler-Maruyama scheme. 
  Here $\sigma(x) \in  \R^{n\times k},$ $\ k\le n$  and $\Sigma(x) =\sigma(x)\sigma^{tr}(x)$ 
  is   non-singular and  positive definite.  
  Note also that the invariant measure is $ \mu(dx) =  Z^{-1} \exp\{-U(x)\}dx$, where $Z$ is the  normalizing constant.
This example demonstrates the applicability of our approach for stochastic dynamics with multiplicative noise. Moreover,   it proves once more that the optimal parameter set $\theta_h^*$ derived from the discrete time analogue of \VIZ{ovdLang}
does not depend on the  discretization time step $h$ as $h\to 0$.

The reason for studying the optimization through the discrete schemes  is
(i) the numerical parametrization and minimization
are done  always in a context of discretized dynamics, 
(ii) it demonstrates clearly the passage from 
Markov chain approximation to continuous time process in the variational problem that defines the best-fit procedure. 
 Indeed, as proved in Theorem~\ref{prop1Lang}, \ref{numerics-sec5},
 when the discretization time step  of the numerical scheme tends to zero the 
 local minima of RER  agree   
with  the results for the continuous time case, Theorem~\ref{REminThm}~b).

%
%
%
%
\subsection{Examples of coarse graining maps}
\label{Center}

\noindent
{\bf Center of mass.}
We consider the linear transformation $\COP_q$, eq. \VIZ{linearMap}, with 
\begin{equation*}
\zeta_{ji} = \BARIT m^{-1}_j m_i  \chi_{C_j}(i)\quad j=1,\dots,M,  \ i=1,\dots,N    \COMMA
\end{equation*}
that maps  $N/M $ microscopic particles to its center of mass (the CG particle), where $m_i$  is the mass of the $i$-th particle, $ \BARIT m_j =\sum_{i\in C_j} m_i$ and  $\chi_{C_j} $ is the indicator function of the set
$C_j=\{i: \text{ particle $i$ contributes to CG particle $j$}\}$.
In order to  simplify the demonstration we consider the special case 
 where $M=2$ and $N$ is even and 
$$
\COP_q =  
\begin{bmatrix}
\BARIT m_1^{-1} m_1 I_3 & \dots & \BARIT m_1^{-1}m_{N/2} I_3& 0 &\dots &0  \\
0 & \dots & 0& \BARIT m_2^{-1}m_{(N/2)+1} I_3& \dots &  \BARIT m_2^{-1}m_{N}I_3
\end{bmatrix}   
$$
where $ \BARIT m_1 = \sum_{i=1}^{N/2} m_i$,  $ \BARIT m_2 = \sum_{i=(N/2)+1}^Nm_i$, and note 
 that $ \text{rank}(\COP_q) = 6 =3 M$, 
and
$$
\COP_p=\cgmass\COP_q \mass^{-1} =
\begin{bmatrix}
   I_3 & \dots &   I_3& 0 &\dots &0  \\
0 & \dots & 0&   I_3& \dots &    I_3
\end{bmatrix} \COMMA
$$
with $\GCOP_p = \mass \GCOP_q \cgmass^{-1} = \COP_p^{tr}$. 
 Thus the weighted norm $\| \cdot \|_{\GCOP_p\Xi} $ appearing in the optimization problems \VIZ{Opt1}-\VIZ{Opt4},
  for the case of identical and constant covariance coefficients $ \sigma$ for which $\Xi = \sigma^{-1} I_{3N}$,  is
  $$ \| u \|_{\GCOP_p\Xi}^2 = \sigma^{-2}  u^{tr} (\GCOP_p )^{tr}  \GCOP_p u = \frac{N}{2\sigma^2}  u^{tr} u\, , \      
 u=(u_1,u_2)^{tr}\in \R^{6} \COMMA $$
    as $\GCOP_p u =(u_1,\dots , u_1, u_2,\dots, u_2)^{tr}$.  The minimization problem \VIZ{Opt1}  becomes
    \begin{equation*}
 \theta^* = \text{argmin}_\theta \EXPECTWRT{\mu}{ \frac{N}{4\sigma^2}
    \sum_{j=1}^2\Big\|   \sum_{i\in C_j}  F_i(q) -\BARIT F_j(\COP_q q;\theta)  \Big\|^2}\COMMA
 \end{equation*}
for $ \BARIT F(\bq;\theta)=(\BARIT F_1(\bq;\theta), \BARIT F_2(\bq;\theta))^{tr}\in \R^{6} $, where $\| \cdot\| $ denotes the Euclidean norm in $\R^3$.

\bigskip
\noindent
{\bf An orthogonal projection.}
A   simple example  of an orthogonal mapping is 
 the projection on the first $M$ coordinates, i.e., for $
 q=(q_1,\dots, q_N)^{tr}$,  $ \COP_q  q = (q_1,\dots,q_M)^{tr}$, represented by  the $3M\times 3N$ matrix
 $$
 \COP_q  = 
\begin{bmatrix}
I_{3M} & 0_{3M\times(3N-3M)}
\end{bmatrix}\COMMA
$$
  and
 $$
\COP_p=\cgmass\COP_q \mass^{-1} =
\begin{bmatrix}
 I_{3M} &  0_{3M\times(3N-3M)}
\end{bmatrix} =\COP_q\COMMA
$$
 with $ \cgmass = \DIAG(m_1I_3, \dots, m_M I_3) $, for which holds
 $\GCOP_p = \COP_p^{tr} $.
When the diffusion  coefficient  $ \sigma$ is  a constant same for all particles  
the weighted norm appearing in the optimization problems \VIZ{Opt1}-\VIZ{Opt2} is
 $$
 \| u \| _{\GCOP_p\Xi} ^2   =\frac{1}{\sigma^2} u^{tr}u, \ \text{ for any }   u \in \R^{3M}  \COMMA
$$
and the minimization problem   
\begin{equation*}
 \theta^* = \text{argmin}_\theta \EXPECTWRT{\mu}{ \frac{1}{2\sigma^2}
   \sum_{i=1}^M  \left\|  F_i(q) -\BARIT F_i(\COP_q q;\theta)  \right\|^2}\PERIOD
 \end{equation*}

\section{Data driven coarse-graining  and path-space variational inference}\label{data_driven}
In the previous sections we built optimal coarse-grained dynamics such as (\ref{CGLang}) as  approximations  of  microscopic stochastic dynamics
in the same general class of stochastic differential equations,
e.g.,   the microscopic  Langevin dynamics (\ref{Lang}). The path-space information approach proposed there, 
allowed us to systematically compare microscopic and coarse-grained  dynamics using the path-space relative entropy approach in Theorem~\ref{REminThm} and set up the corresponding parameter path-space variational inference  problems. 

In this section we explore a different but asymptotically equivalent (in the number of data) perspective, which is {\em data-driven} in the sense  that it treats any available  microscopic simulator 
purely as a means of producing statistical data in the form of time series.

The  primary  new elements of this  coarse-graining approach  are: 
(a) We derive  the  parametrization of the coarse-grained models  $\bar Q^\theta_{0, T]}$ by  optimizing their information 
content (in  the path space) with respect  to available  time series data $\cal D$ coming from 
a fine-scale simulation. We use  computable formulas similar to those 
for the relative entropy rate (RER)  discussed earlier;
(b) In (a) we  do not need that the microscopic scale  time series data $\cal D$ are obtained from 
a model $P_{[0, T]}$ in the same mathematical class as the coarse-grained models or that it is even Markovian.
Thus we do not require a microscopic Langevin model $P_{[0, T]}$ as was done in Section~\ref{RERminS} and  Section~\ref{Langevin-sec}  
or any other  explicitly known  molecular dynamics, at least provided sufficient data is available;
(c)	 
Due to  information inequalities such as \VIZ{CKP} and \VIZ{DKP-bound}, 
the present data-driven, path information-based coarse-graining methodology implies transferability of the parameterization to   different  observables, at a specific thermodynamic point,  by only training  a single observable, namely the relative entropy. 
 Numerical tests for the verification of this observation  is the subject of ongoing work.  We also refer to some recent prior work that relates observables and relative entropy, \cite{DKPP, Tsourtis:2015}. 
(d) The relative entropy rate (RER) approach in (\ref{VP-path2}) which is also reflected in our data-driven approach, see for instance eq. (\ref{Path-L6}) below, allows us to train models to be predictive  at long-time regimes, and not just for any fixed finite time window $[0, T]$.

One of the key points in this method is taking advantage of the ordering 
of the distributions in relative entropy, that allows us to write the path-space variational inference problem as an average of the available fine-scale time series data. In this sense the method relies specifically on the availability of ``big data'' from the microscopic solver. Furthermore, it  provides a systematic approach  to compress them in the form of the coarse-grained model $\bar Q^\theta_{[0, T]}$ with controlled loss of information measured by  RER, in analogy to  information criteria, \cite{Akaike}, \cite{Takeuchi}, see  formula (\ref{VP-CG2}).

For simplicity in the presentation, we focus on discrete in time parametrized coarse dynamics, for instance numerical schemes for Langevin equations, such as the ones considered in
\ref{numerics-sec5}. In analogy to the parametrized dynamics (\ref{CGLang}), we consider 
 the parametrized coarse-grained  transition probabilities $ \BARIT{p}^\theta(\COP X, \COP X')$, i.e., the probability for the CG state $\COP X $ given that the system is at $\COP X$,  that correspond to a discretization scheme for (\ref{CGLang}), see also 
 \ref{numerics-sec5}. 
Then   the corresponding  coarse-grained path-distribution  
is, assuming Markovianity, 
\begin{equation}\label{path-measure-CG}
\bar Q^\theta_{[0, T]}=\bar Q^\theta(\COP X_0,\dots,\COP X_T) = \mu(\COP X_0) \BARIT{p}^\theta(\COP X_0, \COP X_1)\dots \BARIT{p}^\theta(\COP X_{T-1}, \COP X_T)\, ,
\end{equation}
where $\mu$ denotes the initial distribution of the process and $\BARIT{\mathcal{D}}=\{\COP X_1, \COP X_2,...,\COP X_T \}$ is a typical coarse-grained time series corresponding to the microscopic time series $\mathcal{D}=\{X_1, X_2,...,X_T \}$.
 
\medskip
\noindent{\em Finite-time regime:}
We first  consider an ensemble of   fine-scale  data given in the form of $M$  time series   
$\mathcal{D}^k=\{X^k_1, X^k_2,...,X^k_T \}$, $k=1, ..., M$,  obtained  from a fine-scale 
molecular simulation algorithm, up to a prescribed time horizon $T$. 
The corresponding coarse space time series is $\BARIT{\mathcal{D}}^k=\{\COP X^k_1, \COP X^k_2,...,\COP X^k_T \}$,
$k=1, ..., M$. We note here that in this case we do not necessarily need that this data set is obtained 
from a Langevin or any other Markovian algorithm.
Then the {\em typically unknown}  path-space distribution of this coarse-space time series  is denoted  by 
$$\BARIT{P}=\BARIT{P}(\COP X_0,\dots,\COP X_T)\PERIOD$$ 
Note that $\BARIT{P}$ is the push forward of the microscopic measure $P$, see relation \VIZ{defPFW} Section~\ref{InfoCG}.
Furthermore $\BARIT{P}$ is computationally accessible from the microscopic simulation which samples $P$ and projects the $X_i$'s onto $\COP X_i$'s.

On the other hand we also consider the parametrized coarse-grained family,  defined in (\ref{path-measure-CG}),
$ \bar Q^\theta =\bar Q^\theta(\COP X_0,\dots,\COP X_T)$. In order to obtain the optimal parametrized coarse-grained transition probabilities  $\BARIT{p}^{\theta^*}(\COP X, \COP X')$, 
for a parameter vector $\theta=\theta^*$,  we need to minimize the path-space relative entropy,
i.e., consider the minimization problem 
\begin{equation}
\label{min-RE}
\theta_T^*=\arg\min_\theta \mathcal{R}(\BARIT{P}|\bar Q^\theta)\PERIOD
\end{equation}
 Furthermore, for the path-space relative entropy we have by the law of large numbers that
\begin{equation*}
\label{MLE0}
\RELENT{\BARIT{P}}{\bar Q^{\theta}}=\lim_{M \to \infty} \hat{\mathcal R}_M(\BARIT{P} |\, \bar Q^{\theta})
\end{equation*}
where we define the unbiased estimator for the relative entropy   
\begin{equation}\label{Path-Lestim}
\hat{\mathcal R}_M(\BARIT{P} |\, \bar Q^{\theta}):=\frac{1}{M}\sum_{k=1}^M\log\frac{\BARIT{P}(\COP X^k_1, \COP X^k_2,...,\COP X^k_T)}{ 
	\bar Q^\theta({\COP X^k_1, \COP X^k_2,...,\COP X^k_T})}\, .
\end{equation}
Therefore, the minimization principle (\ref{min-RE})
 becomes
\begin{equation}\label{MLE}
\begin{aligned}
\arg\min_\theta \hat{\mathcal R}_M(\BARIT{P} |\, \bar Q^{\theta})= & \arg\max_\theta \frac{1}{M}\sum_{k=1}^M\log 
\bar Q^\theta({\COP X^k_1, \COP X^k_2,...,\COP X^k_T})
 \\
& - \frac{1}{M}\sum_{k=1}^M\log \BARIT{P}(\COP X^k_1, \COP X^k_2,...,\COP X^k_T) \COMMA
\end{aligned}
\end{equation}
which does not require   a priori the knowledge of the microscopic probability distribution 
or its push-forward $\BARIT{P}(\COP X^k_1, \COP X^k_2,...,\COP X^k_T)$.
Therefore, we obtain  from (\ref{MLE}) 
and (\ref{path-measure-CG}) the following   maximization principle 
\begin{equation}\label{Path-L}
\theta_T^* \approx \arg\max_\theta\sum_{i=1}^{T}\sum_{k=1}^M\log 
\BARIT{p}^\theta(\COP X^k_i,\COP X^k_{i+1})
 + \sum_{k=1}^M\log\mu(\COP X^k_0)
\PERIOD
\end{equation}
For a time window $T\gg 1$ the last term in (\ref{Path-L}) that involves the initial data becomes  negligible, therefore we obtain a coarse-grained path-space likelihood maximization principle and the corresponding estimator 
$\hat \theta(M, T)$ of $\theta^*$:
\begin{equation}\label{Path-L2}
\hat \theta(M, T)=\arg\max_\theta
L(\theta; \{\COP X^k_i\}_{i=1, k=1}^{T, M})\, ,
\end{equation}
where 
\begin{equation}\label{Path-L3}
L(\theta; \{\COP X^k_i\}_{i=1, k=1}^{T, M})=\sum_{i=1}^{T}\sum_{k=1}^M\log 
\BARIT{p}^\theta(\COP X^k_i,\COP X^k_{i+1})
\PERIOD
\end{equation}

Note that if the transition probabilities in \VIZ{Path-L3}  are replaced with a stationary 
measure and $T$ corresponding to
independent samples $\mathcal{D}=\{X_1, X_2,...,X_T \}$,
then \VIZ{Path-L2} becomes the classical Maximum Likelihood Principle (MLE), \cite{CaseBerg:01}.
In this sense \VIZ{Path-L} is a  maximum likelihood for the  coarse-grained  time series,    
$\BARIT{\mathcal{D}}=\{\COP X_1, \COP X_2,...,\COP X_N \}$
of the fine-scale  process, and thus includes dynamics information and temporal correlations.

\medskip
\noindent{\em Stationary regime:}
If the time series associated with the data set $\mathcal{D}^k=\{X^k_1, X^k_2,...,X^k_T \}$, $k=1, ..., M$ are stationary, then they are statistically indistinguishable and we will eventually drop the index $k$, referring to the data set as 
 $\mathcal{D}=\{X_1, X_2,...,X_T \}$.
In this case  the estimator in (\ref{Path-L3}) simplifies significantly. Indeed, for $M\gg 1$, we have that 
\begin{equation}\label{Path-L4}
L(\theta; \{\COP X^k_i\}_{i=1, k=1}^{T, M})=\sum_{i=1}^{T}\sum_{k=1}^M\log 
\BARIT{p}^\theta(\COP X^k_i,\COP X^k_{i+1}) \approx M \sum_{i=1}^{T} \E_{\BARIT P}\big[\BARIT{p}^\theta(\COP X_i,\COP X_{i+1})\big]
\, ,
\end{equation}
  $\E_{\BARIT{P}}$ denotes the expectation with respect to the path   distribution  
  $\BARIT{P}(\COP X_1, \COP X_2,...,\COP X_T)$. 
Using the stationarity of the time series in  (\ref{Path-L4}) we have that 
\begin{equation*}\label{Path-L5}
\sum_{i=1}^{T} \E_{\BARIT{P}}\big[\BARIT{p}^\theta(\COP X_i,\COP X_{i+1})\big]= 
T \E_{\BARIT P}\big[\BARIT{p}^\theta(\COP X_1,\COP X_{2})\big]
\PERIOD
\end{equation*}
However, due to stationarity we have the unbiased estimator for 
$\E_{\BARIT P}\big[\BARIT{p}^\theta(\COP X_1,\COP X_{2})\big]$, \cite{KP2013}, \cite{PK2013},
using a single time series $\mathcal{D}=\{X_1, X_2,...,X_T \}$:
$$
\E_{\BARIT P}\big[\BARIT{p}^\theta(\COP X_1,\COP X_{2})\big] \approx \frac{1}{T} \sum_{i=1}^{T} \BARIT{p}^\theta(\COP X_i,\COP X_{i+1}) \, .
$$
Therefore the stationary analogue of (\ref{Path-L2}) is
\begin{equation}\label{Path-L6}
\theta^* \approx \hat \theta_s(T)=\arg\max_\theta
L_s(\theta; \{\COP X_i\}_{i=1}^{T})\, ,
\end{equation}
where for the stationary time series $\mathcal{D}=\{X_1, X_2,...,X_T \}$ we define 
\begin{equation}\label{Path-L7}
L_s(\theta; \{\COP X_i\}_{i=1}^{T})=\sum_{i=1}^{T}\log 
\BARIT{p}^\theta(\COP X_i,\COP X_{i+1})
\PERIOD
\end{equation}

\medskip
\noindent{\em Time dependent regime:}
	The optimization principle in (\ref{Path-L2}) can be also easily extended so that we can obtain a further  improved but computationally more costly  time-dependent optimal parametrization for (\ref{min-RE}) where now we seek
	$$
	\theta^*=\theta^*(i)\, , \quad \mbox{where} \quad  0\le i \le T\, ,
	$$
	solving the approximate optimization problem for $\theta=\theta(i)$, $i=0,..., T$:
	$$
	\arg\min_{\theta}\sum_{i=0}^{T-1}\sum_{k=1}^M\log 
	\BARIT{p}^{\theta(i)} (\COP X^k_i,\COP X^k_{i+1})\, .
	$$
	In \cite{Kalliadasis2015} the authors also obtained data-driven parametrization of multi-scale diffusions on a finite time window $[0, T]$ based on minimizing 
	(\ref{IMC1}) for particular observables $\phi_i$. On the other hand, in \cite{KP2013} a similar  
	problem was considered 
	for  multi-scale diffusions, where coarse-grained parametrization was  based on minimizing 
	instead information metrics, e.g.,  \VIZ{VP-path2}.

\begin{remark}
{\rm 
	Confidence intervals for the estimator $\hat \theta_s(T)$ in (\ref{Path-L7})
	can be provided in terms of the asymptotic normality (in  $T\gg 1$)  of the estimator.
	The corresponding  (Gaussian) asymptotic variance is given in terms of the inverse of 
	the path-space Fisher information matrix (FIM). We refer to
	\cite{KP2013,PVK2013,Tsourtis:2015} for a discussion of such results, as well as the definition of the path FIM and its use for sensitivity analysis  of  stochastic dynamics that include Langevin dynamics and kinetic Monte Carlo algorithms.
	}
\end{remark}

\section{Conclusions}\label{discussion}
In this work we presented a thorough examination of coarse-graining of {\it non-equilibrium} molecular 
systems using {\it path-wise  information metrics}. We have introduced the minimization problem for 
optimizing coarse models based on  relative entropy  for comparing continuous time diffusion processes.
The derived scheme is similar to the widely applied  force-matching  method used in computational 
coarse-graining
which, however, is restricted to equilibrium processes.

The main novelties of the proposed approach are summarized in the following points:
(a) It is applicable to transient regimes of non-equilibrium processes, since it directly 
involves information along the whole path-space;  
(b) It connects the path-space relative entropy minimization with an (extended) force matching  problem for continuous time dynamics.
(c) It becomes entirely data-driven when the microscopic dynamics are replaced with corresponding correlated data in the form of time series.
From a more general perspective the proposed scheme is directly related to dynamical data fitting as well as to machine learning algorithms. Indeed, the path-space information approach allowed us to relate the 
RER minimization problem to corresponding parameter optimization problems, obtained from data-driven 
methodologies, in the sense that it treats the microscopic simulator as means of producing statistical data in the form of time-series;
(d) The interpretation of the dynamics with continuous time process  demonstrates that the  RER for stationary dynamics is independent of the time step  for any numerical discretization scheme.  Most importantly the RER perspective shows that the corresponding optimization method is extendable to infinite times and non-reversible systems as is demonstrated in the current study for continuous time diffusion processes and in \cite{KP2013} for Markov chains;
(e)The approach is generally applicable to stochastic dynamics such as  Kinetic Monte Carlo algorithms and reaction networks.
 
Current work concerns the numerical application of the proposed RER minimization methodology to coarse-graining of molecular systems under equilibrium and non-equilibrium conditions.

\section*{Acknowledgements}
The research of E.K. and V.H. was supported by the European Union (ESF) and Greek national funds through the Operational Program Education and Lifelong Learning of the NSRF-Research Funding Program: THALIS. 
The research of  V.H. was also partially supported by ARISTEIA II.
The research of M.K. was supported in part by the Office of Advanced Scientific Computing Research, U.S. Department of Energy under Contract No. DE- SC0010723. The research of P.P. was partially supported by the U.S. Department of Energy Office of Science, Office of Advanced Scientific Computing Research, Applied Mathematics program under Award Number DE-SC-0007046. 

\appendix

 \section{Proofs}\label{appendixProof}
\subsection{Proof of Theorem~\ref{PropRecon}}\label{appendixProof1}
We consider a coarse  observable function $f(x)=g(\COP x)$,  where $f\in C^2(\R^n;\R)$ and $g\in C^2(\R^m;\R)$, and such that the conditions of Theorem~\ref{PropRecon} are guaranteed.  $ C^2(\R^n;\R)$ denotes the space of twice differentiable functions $f: \R^n\to \R$.

From the martingale problem,   \cite[Section 8.3]{oksendal2003}, \cite{karatzas1991},
 given $\tb, \tilde\Sigma$, defined in relations \VIZ{backdrift} and \VIZ{backdiff},  there exists unique process $\{\bX_t\}_{t>0}$   which is the  solution of \VIZ{rSDE1}  such that  
 $$\bX_t  =  \COP \tX_t, \quad \text{ in distribution.}$$   
Denote
$$\Ll f(x) = \sum_{i=1}^n \tb_i(x) \frac{\partial f}{\partial x_i}  (x) +  \frac12 \sum_{i,j}^n  \tilde\Sigma_{ij}(x;\theta) \frac{\partial^2 f}{\partial x_i \partial x_j}(x)    $$
  with generator of $X_t$.
For $f(x)=g(\COP x)$, using relations \VIZ{backdrift} and \VIZ{backdiff},  see also \cite[Lemma 7.3.2]{oksendal2003},   we can write
$$\Ll f(x) =   \sum_{I=1}^m (\bbr (\COP x))_I\frac{\partial g}{\partial \bx_I}  (\COP x) + \frac12 \sum_{I,J=1}^m  \BARIT \Sigma_{IJ}(\COP x;\theta)\frac{\partial^2 g}{\partial \bx_I \partial \bx_J}     (\COP x) ds$$
Let $\tau>0$ be a stopping time,
from  Ito's formula 
\begin{eqnarray*}
f(\tX_t) &=& f(x) + \int_0^t \Ll f(\tX_s) ds + \int_0^t \nabla f(\tX_s)^{tr} \tilde \sigma dB_s 
\end{eqnarray*}
where  $X_0 = x$,
 by applying expectation on both sides, 
\begin{eqnarray*} 
\E^x[f(\tX_\tau)] &=& f(x) +  \E^x\left[ \int_0^t \Ll f(\tX_s) ds + \int_0^t \nabla f(\tX_s)^{tr} \tilde \sigma dB_s  \right]\\ 
 &=&  g(\COP x) +   \E^{\COP x}\left[\int_0^\tau \sum_{I=1}^m (\bbr (\bX_s))_I\frac{\partial g}{\partial \bx_I}  (\bX_s)  \right.\\
 &&\left. \qquad \quad +\frac12 \sum_{I,J=1}^m  \BARIT \Sigma_{IJ} (\bX_s;\theta)\frac{\partial^2 g}{\partial \bx_I \partial \bx_J}     (\bX_s) ds \right] 
\end{eqnarray*}

Thus  
$$ \E^x[f(\tX_\tau)] = \E^{\COP x}[g(\bX_\tau)], $$
for any coarse observable $f(x)$.
 \hfill $\square$  
\subsection{Girsanov Theorem}\label{appendixProof2}
  The Girsanov theorem states the conditions under which the path space measures 
  $\PATHPT$ and $\tPATHPT$
  are absolutely continuous and provides a closed form of  the Radon-Nikodym density $\frac{d\PATHPT}{d\tPATHPT}(X_t,t)$.  
Suppose that there exists a process $\{u(X_s;\theta)\}_{s\ge 0} $   in $\R^{k}$ such that
\begin{equation*} 
 \sigma(X_s)u(X_s;\theta) = b(X_s) -  \tb(X_s;\theta)\COMMA
 \end{equation*}
and satisfies Novikov's condition $\E\left[ \EXP{\frac12\int_0^T |u(X_s;\theta)|^2\,ds}\right] <\infty $, \cite{oksendal2003}. We define
\begin{eqnarray*} 
M_t:=\frac{d\mu_0}{d\nu_0}(X_0)\exp\left\{-\int_0^t \SCPROD{u(X_s;\theta)}{dB_s} -\frac12 \int_0^t |u(X_s;\theta)|^2\,ds\right\}\COMMA
\end{eqnarray*} 
where $ \SCPROD{u(X_s;\theta)}{dB_s}=\sum_{i=1}^k  u_i(X_s;\theta)  dB^i_s$ and 
$|u(X_s;\theta)|^2 =\sum_{i=1}^k  u^2_i(X_s;\theta)$.
Then the Girsanov theorem  yields that 
 $\tPATHPT$ is   absolutely continuous with respect to $ \PATHPT$,  $\tPATHPT~\!\ll~\!\PATHPT$, and 
\begin{eqnarray}\label{RN}
\frac{d\PATHPT}{d\tPATHPT}(X_t,t) = M_t\PERIOD
\end{eqnarray} 
Furthermore,  the process $\hat B_t:= \int_0^t u(X_s;\theta)\,ds + B_t $ is a $k$-dimensional Brownian motion with respect to $\PATHPT$.

\subsection{Proof of Theorem~\ref{thmRel}}\label{appendixProof3}
\begin{enumerate}[a)]
\item
The Novikov condition  $$\E\left[ \EXP{\frac12\int_0^T |u(X_s;\theta)|^2 ds}\right] <\infty$$ ensures that
 $\RELENT{\PATHPT}{\tPATHPT} < \infty $. 
Thus 
\begin{eqnarray*}
 &&\RELENT{\PATHPT}{\tPATHPT} - \Hh^T(\PATHPT\SEP\tPATHPT)  - \RELENT{\mu_0}{\nu_0} = \\
 &&= \EXPECTWRT{\PATHPT}{\log\frac{d\PATHPT}{d\tPATHPT}} - \EXPECTWRT{\PATHPT}{\frac12 \int_0^T |u(X_s;\theta)|^2 ds} - \RELENT{\mu_0}{\nu_0}\\
 &&= \EXPECTWRT{\PATHPT}{-\int_0^T \SCPROD{u(X_s;\theta)}{dB_s}  - \int_0^T |u(X_s;\theta)|^2\, ds} \\
 &&= \EXPECTWRT{\PATHPT}{-\int_0^T \SCPROD{u(X_s;\theta)}{dB_s  + u(X_s;\theta)\,ds}} = 
     \EXPECTWRT{\PATHPT}{-\int_0^T \SCPROD{u(X_s;\theta)}{d\hat B_s}} =0\COMMA
\end{eqnarray*}
the last term is zero  since $d\hat B_s$ is a Brownian motion with respect to $\PATHPT$.

\item
Since $\PXT$ is stationary 
$$
  \Hh^T(\PATHPT \SEP\tPATHPT) =  T \EXPECTWRT{\mu}{\frac12 |u(X;\theta) |^2}  \COMMA
$$
thus, based on the representation \VIZ{RE} of the relative entropy,
\begin{eqnarray*}
 \lim_{ T\to \infty} \frac{1}{T} \RELENT{\PATHPT}{\tPATHPT}
 &=&  \lim_{ T\to \infty} \left[ \frac{1}{T} \Hh^T(\PATHPT\SEP\tPATHPT) + \RELENT{\mu}{\nu}\right]\\
 &=&  \EXPECTWRT{\mu}{\frac12 |u(X;\theta)|^2}  \PERIOD
\end{eqnarray*}
Recalling  the definition of the RER, 
$$
 \ENTRATE{P}{Q^{\theta}}  = \lim_{ T\to \infty} \frac{1}{T} \RELENT{\PATHPT}{\tPATHPT} \COMMA
$$
we conclude that 
$\ENTRATE{P}{Q^{\theta}}   =   \EXPECTWRT{\mu}{\frac12 |u(X;\theta)|^2}$.
\end{enumerate}

\subsection{Proof of Corollary~\ref{thm3}}\label{appendixProof4}
 
By assumption there exists $  u(x;\theta) $ such that 
\begin{eqnarray*}
           \sigma(x)u(x;\theta) = b(x) -  \tb(x;\theta), \text{ for all } x\in\R^n\PERIOD
\end{eqnarray*}
Since $\text{rank}(\sigma(x))=r=k$, i.e., $\sigma(x)$ has full rank, then $ \sigma^{tr}(x) \sigma(x) \in \R^{k\times k}$ is invertible, thus    a  solution $ u(x;\theta)\in \R^k $  is given by
\begin{eqnarray*}
           u(x;\theta) = \left[\sigma^{tr}(x)\sigma(x)\right]^{-1}\sigma^{tr}(x)\left(b(x) -  \tb(x;\theta)\right) = \Xi(x) \left(b(x) -  \tb(x;\theta)\right) \COMMA
\end{eqnarray*}
where $\Xi(x) = \left[\sigma^{tr}(x)\sigma(x)\right]^{-1}\sigma^{tr}(x)$
and 
\begin{eqnarray*}
          | u(x;\theta)|^2&=& u^{tr}(x;\theta)u(x;\theta) \\
           &=&\left(b(x) -  \tb(x;\theta)\right)^{tr} \Xi^{tr}(x)\Xi(x) \left(b(x) -  \tb(x;\theta)\right) = \| b(x) -  \tb(x;\theta)\|_{\Xi}^2\COMMA
\end{eqnarray*}
From Theorem~\ref{thmRel}~b)    the RER is  
   $$
   \ENTRATE{\PATHP}{\tPATHP} =\EXPECTWRT{\mu}{\frac12 |u(X;\theta)|^2}\COMMA
   $$
  thus   substituting of the previously derived   form of $u(x;\theta)$ we prove~\VIZ{RERf1},
    $$
   \ENTRATE{\PATHP}{\tPATHP} =\EXPECTWRT{\mu}{\frac12\| b(x) -  \tb(x;\theta)\|_{\Xi}^2}\COMMA
   $$ 
   where the norm $\| \cdot \|_{\Xi}$ is defined in \VIZ{XiNorm}.
 
   \subsection{Proof of Theorem~\ref{REminThm}}\label{appendixProof5}
\begin{enumerate}[a)]
\item
From Corollary~\ref{thm3}(\ref{FTRE}),  we have that 
\begin{equation*} 
\RELENT{\PATHPT}{\tPATHPT} =\Hh^T(\PATHPT\SEP\tPATHPT)  + \RELENT{\mu_0}{\nu_0}\COMMA
\end{equation*} 
where
\begin{equation*} 
 \Hh^T(\PATHPT\SEP\tPATHPT) =  
 \EXPECTWRT{\PATHPT}{\frac12 \int_0^T  \|  b(X_s) -  \tb(X_s;\theta)\|^2_{\Xi}   \,ds}  \PERIOD
\end{equation*}
As $\mu_0, \nu_0$,  are independent of $\theta$,  we only need to  prove that
\begin{equation*}
 \nabla_{\theta}   \EXPECTWRT{\PATHPT}{\frac12 \int_0^T\! \|  b(X_s) -  \tb(X_s;\theta)\|^2_{\Xi}   \,ds}\!=\! \nabla_{\theta}  \EXPECTWRT{\PATHPT } { \frac12 \int_0^T\!  \| \COP b(X_s) - \bbr(\COP X_s;\theta) \|_{\GCOP\Xi}^2\,ds }\! .
 \end{equation*}
Recall the representation \VIZ{brecRep} 
$$ 
 \tb(x;\theta) =\GCOP\bbr(\COP x;\theta) + (I-\GCOP\COP)y^{\perp}(x)\COMMA\;\;\;\mbox{and}\;\;
 \COP(I-\GCOP\COP)y^{\perp}(x)  = 0\COMMA$$
for all  $x\in \R^n$. We also write 
$  b(x) =\GCOP \COP b(x) + (I-\GCOP\COP)b(x) $ 
and  have 
\begin{eqnarray*}
 \| b(x) - \tb(x;\theta) \|_\Xi^2& = & \|\GCOP\COP \left(b(x) - \tb(x;\theta)\right) \|_\Xi^2 + \|(I-\GCOP\COP )\left(b(x) - y^{\perp}(x) )\right) \|_{\Xi}^2 \\
 & = &
 \|\COP \left(b(x) - \tb(x;\theta)\right) \|_{\GCOP\Xi}^2  +\|(I-\GCOP\COP )\left(b(x) - y^{\perp}(x) )\right) \|_{\Xi}^2\COMMA
\end{eqnarray*}
for all  $x\in \R^n$, where the cross terms are zero as
$\COP(I-\GCOP\COP) y^{\perp}(x)  =0$.
Thus from the assumption that 
$ y^{\perp}(x)  $  is independent of $\theta$ we have
\begin{eqnarray*}
\nabla_\theta \RELENT{\PATHPT}{\tPATHPT} &=& \nabla_\theta \Hh^T(\PATHPT\SEP\tPATHPT)\\
&=& \nabla_\theta \EXPECTWRT{\PATHPT } { \frac12 \int_0^T \| \COP b(X_s) - \bbr(\COP X_s;\theta) \|_{\GCOP\Xi}^2\,ds } \COMMA
\end{eqnarray*}
and
 \begin{equation*}
\mathrm{argmin}_{\theta\in\Theta}\RELENT{\PATHPT}{\tPATHPT} = \mathrm{argmin}_{\theta\in\Theta}  
  \EXPECTWRT{\PATHPT } { \frac12 \int_0^T \| \COP b(X_s) - \bbr(\COP X_s;\theta) \|_{\GCOP\Xi}^2\,ds } \PERIOD
 \end{equation*}
\item 
We recall that  for stationary process $\PXT $, see Corollary~\ref{thm3}
$$\ENTRATE{\PATHP}{\tPATHP} =\frac{1}{T}\Hh^T(\PATHPT\SEP\tPATHPT)  \COMMA
$$
  and from a)  we have
$$ 
 \nabla_\theta \Hh^T(\PATHPT\SEP\tPATHPT)\\
= \nabla_\theta \EXPECTWRT{\PATHPT } { \frac12 \int_0^T \| \COP b(X_s) - \bbr(\COP X_s;\theta) \|_{\GCOP\Xi}^2\,ds }\COMMA
 $$
thus
\begin{eqnarray*} 
   \nabla_\theta \ENTRATE{\PATHP}{\tPATHP} &=& \frac{1}{T}\nabla_\theta \EXPECTWRT{\PATHPT } { \frac12 \int_0^T \| \COP b(X_s) - \bbr(\COP X_s;\theta) \|_{\GCOP\Xi}^2\,ds } \\
   &=&   \frac{1}{T}T\nabla_\theta \EXPECTWRT{\mu} { \frac12 \| \COP b(X ) - \bbr(\COP X;\theta) \|_{\GCOP\Xi}^2  }\\
 &=&  \nabla_\theta \EXPECTWRT{\mu} { \frac12  \| \COP b(X ) - \bbr(\COP X;\theta) \|_{\GCOP\Xi}^2  }\COMMA
   \end{eqnarray*}
   where we used that $\PXT$ is stationary with the invariant measure $\mu(dx)$.
Thus
 \begin{equation*}
\mathrm{argmin}_{\theta\in\Theta}\ENTRATE{\PATHP}{\tPATHP} =  
\mathrm{argmin}_{\theta\in\Theta}\EXPECTWRT{\mu} {\frac12  \| \COP b(X) - \bbr(\COP X;\theta) \|_{\GCOP\Xi}^2 }\PERIOD   \qquad \qquad \square 
  \end{equation*}  

\end{enumerate}

\section{Relative entropy rate minimization for numerical schemes}\label{numerics-sec5}

The   relative entropy rate  $\ENTRATE{\PATHP}{\tPATHP}$, for discrete time Markov chains  $\{x_i\}_i$,$\{\tilde  x_i\}_i$
with transition probabilities $p^h(x,x')$ and $q^h_\theta(x,x')$ respectively is,  
  \cite{KP2013},  
\begin{equation}\label{discrRER}
\mathcal H(\PATHP|\tPATHP) =\frac1h  \int\int
 \mu(x) p^h(x,x')\log \frac{p^h(x,x')}{q^h_\theta(x,x')}dx'dx\PERIOD
\end{equation}
The RER is  described here for the specific (stationary) Markov chains 
generated by the numerical schemes for the Langevin \VIZ{Lang} and the overdamped Langevin \VIZ{ovdLang}
thus formula \VIZ{discrRER} 
is a statistical estimator of  the RER, see  \cite{KP2013}. 

In terms of   the description in Section~\ref{data_driven},  where the discrete time version of the RE is introduced on the coarse space for $\CGPATHP$, eq. \VIZ{Path-Lestim},
here we use the same description though on the  microscopic space
  where $\tPATHP$ is a reconstruction  of $\CGPATHP$.

\subsection{Euler discretization for over-damped Langevin}
The Euler discretization scheme for \VIZ{ovdLang}  with time step $h$  defines the discrete stochastic system 
\begin{eqnarray*}
x_{i+1} = x_i -\frac12 \Sigma(x_i)\nabla \IP (x_i) h + \frac12 \nabla \Sigma(x_i) h+\sigma(x_i)\Delta W_i\COMMA
\end{eqnarray*}
with the solution given by  the Markov chain $\{x_i\}_{i\geq 0}$ on the state space $\R^n$,
  $\Delta W_i~\sim~N(0,h I_n)$ are normally distributed increments.
The transition probability density for the chain $\{x_i\}_i$ is 
\begin{equation*}
p^h(x,x') = \frac{1}{Z^h(x)}\exp\left\{ -\frac{1}{2h}(x' -\Delta x)^{tr} \Sigma^{-1}(x)(x' -\Delta x)\right\}\COMMA
\end{equation*}
where we denote 
$$
\Delta x = x - \frac12 \Sigma(x) \nabla \IP(x) h + \frac12 \nabla \Sigma(x) h\COMMA
$$ 
and $ Z^h(x) = (2\pi h )^{n/2} |\Sigma(x)|^{1/2}$,
with the notation $|A|$ for the determinant of the matrix $A$.

Consider the linear CG map  $\COP:\R^n\to \R^m $, which for  simplicity we assume it is an orthogonal projection from the state space $\R^n$ onto  $\R^m$ such that 
\begin{equation}\label{}
x = \COP x + (I-\COP) x= \BARIT x + \hat  x\PERIOD
\end{equation}
 Note that we use the same letter $\COP x \in \R^m$ for denoting the representation of   $\COP $ in $\R^n$.
  
The reduced model   based on the CG mapping $\COP$ is given by $ \{\BARIT x_i\}$, approximating  the projected Markov chain  
$\{\COP x_i \}$, satisfying 
\begin{eqnarray*}
\BARIT x_{i+1} =\BARIT x_i -\frac12 \BARIT\Sigma(\BARIT x_i;\theta)\nabla \BARIT \IP (\BARIT x_i;\theta) h + \frac12 \nabla \BARIT\Sigma(\BARIT x_i;\theta) h+\Delta \BARIT W_i\COMMA
\end{eqnarray*}
where $\Delta \BARIT W_i \sim N(0,h \BARIT \Sigma(\BARIT x_i;\theta))$.
Hence the transition probability density for the CG chain $\{\BARIT x_i\}_i$ is 
\begin{equation*}
p^h_\theta(\BARIT x,\BARIT x') = \frac{1}{\BARIT Z^h(\BARIT x;\theta)}\exp\left\{ -\frac{1}{2h}(\BARIT x' -\Delta_\theta \BARIT x)^{tr} \BARIT\Sigma^{-1}(\BARIT x;\theta)(\BARIT x' -\Delta_\theta \BARIT x)\right\}\COMMA
\end{equation*}
where 
$$
\Delta_\theta \BARIT x = \BARIT x_i-\frac12 \BARIT\Sigma(\BARIT x_i;\theta)\nabla \BARIT U(\BARIT x_i;\theta) h 
   +  \frac12 \nabla \BARIT\Sigma(\BARIT x_i;\theta) h\COMMA
$$ 
and $ \BARIT Z^h(\BARIT x;\theta) = (2\pi h )^{m/2} |\BARIT\Sigma(\BARIT x_i;\theta) |^{1/2}$.
We consider the  reconstructed  chain $\{\tilde x_i\}_{i\geq 0}$  of  $\{\BARIT x_i\}_{i\geq 0}$ that has the transition probability density 
$$ 
q^h_\theta(x,x')= \nu(x'|\BARIT x ') p_\theta(\BARIT x,\BARIT x')\COMMA
$$ 
where $\nu$ is a probability measure associated with the reconstruction, that we assume independent of the parameter $\theta$.
Note  that we can write the covariance matrix 
$\Sigma(x)$ as 
$$
\Sigma(x)=  
\begin{bmatrix} 
 \Sigma_{11}(x) & \Sigma_{12}(x)\\ 
  \Sigma_{21}(x) &\Sigma_{22}(x)
\end{bmatrix} 
$$
with $ \Sigma_{11}(x) = \COP \Sigma(x) \COP^{tr}\in \R^{m\times m}$,
$\Sigma_{22}(x) = \COP^\perp \Sigma(x) \COP^{\perp,tr}\in \R^{(n-m)\times(n-m)}$, $\COP^{\perp} = \left(I - \COP\right)$
and $\Sigma_{21}^{tr}(x)  =\Sigma_{12}(x)$.
 Then we can rewrite $p^h(x,x')$ as 
\begin{equation}\label{pFactorize}
p^h(x,x')= p_1(x,\BARIT x') p_2(x,\hat x ' |\BARIT x') \COMMA
\end{equation}
where
\begin{eqnarray*}
&&p_1(x,\BARIT x')  = \frac{1}{Z_1^h(x)} \exp\left\{-\frac{1}{2h} (\BARIT x' - \COP \Delta x)^{tr} 
\Sigma_{11}^{-1}(x)   (\BARIT x' - \COP \Delta x)\right\}\COMMA
\end{eqnarray*}
$Z_1^h(x) =  (2\pi h )^{m/2} |\Sigma_{11}(x)|^{1/2}$ and
\begin{eqnarray*}
&&p_2(x,\hat x ' |\BARIT x')  = \frac{1}{Z_2^h(x)} \exp\left\{-\frac{1}{2h} (\hat x' - g(x,\BARIT x'))^{tr} \Xi^{-1}(x)   (\hat x' - g(x,\BARIT x'))\right\}\COMMA
\end{eqnarray*}
where
\begin{eqnarray*}
&&\Xi(x) = \Sigma_{22}(x) - \Sigma_{12}^{tr}(x)\Sigma_{11}^{-1}(x)\Sigma_{12}(x) \in \R^{(n-m)\times(n-m)}\COMMA\\
 &&g(x,\BARIT x') = \COP^\perp \Delta x + \Sigma_{12}^{tr}(x) \Sigma_{11}^{-1}(x) \left(\BARIT x' - \COP \Delta x \right)\COMMA
\end{eqnarray*}
$\Sigma_{22}(x) = \COP^\perp \Sigma(x) \COP^{\perp,tr}\in \R^{(n-m)\times(n-m)}$,
$\Sigma_{21}^{tr}(x)  =\Sigma_{12}(x)$ 
and $Z_2^h(x) =  (2\pi h )^{(n-m)/2} |\Xi(x)|^{1/2}$.

The variational problem for the best-fit of parameters in terms of the relative entropy rate between
the  Markov chains   $\{x_i\}_{i\geq 0}$   and  $\{\tilde x_i\}_{i\geq 0}$  
is demostrated by the following theorem 
revealing   its relation  to a  weighted force-matching optimization.
\begin{thm}\label{thm2}
\begin{enumerate}[a)] Given $h>0$
\item 
\hspace{1cm }$ \mathrm{argmin}_{\theta\in\Theta} \mathcal{H}(\PATHP|\tPATHP)  = 
\mathrm{argmin}_{\theta\in\Theta} \left[\frac1h A(\theta) + B(\theta)\right]
$,

where
\begin{eqnarray*} \label{AB1}
&&A(\theta)=  \frac12 \int \left[-\log |\COP\Sigma(x)\COP^{tr}\BARIT\Sigma^{-1}(\BARIT x;\theta)|
 + \mathrm{Tr}\left( \COP \Sigma(x)\COP^{tr} \BARIT\Sigma^{-1}(\BARIT x;\theta)\right) \right] \,\mu(x) dx \nonumber \\
&&B(\theta)=  \frac{1}{2}\int \left[\left(\COP b(x)-  \BARIT b(\BARIT x;\theta)\right)^{tr}
 \BARIT\Sigma^{-1}(\BARIT x;\theta)\left(\COP b(x)-  \BARIT b(\BARIT x;\theta)\right)\right]\,\mu(x) dx
\end{eqnarray*}
with  
\begin{eqnarray*}
&& b(x) =  -  \frac12\Sigma(x)\nabla U(x)  + \frac12 \nabla \Sigma(x) \COMMA \\
&& \BARIT b(\BARIT x;\theta)=  -\frac12 \BARIT\Sigma(\BARIT x;\theta)\nabla \BARIT U(\BARIT x;\theta)   +  \frac12 \nabla \BARIT\Sigma(\BARIT x;\theta) \PERIOD
\end{eqnarray*}
\item In the limit $h\to 0$ we have
\begin{equation}\label{argmin2}
 \lim_{h\to 0} \min_{\theta\in\Theta} \left[ \frac1h A(\theta) + B(\theta)\right] = \min_{\theta\in\Theta} \left[E(\theta) + m/2\right]
\end{equation}
where
$ E(\theta)= B(\theta)|_{\BARIT \Sigma = \PrS}$, that is 
\begin{eqnarray*}
&&E(\theta)=  \frac{1}{2}\int\mu(x) \left(\COP b(x)-  \BARIT b(\BARIT x;\theta)\right)^{tr}
 (\PrS)^{-1}\left(\COP b(x)-  \BARIT b(\BARIT x;\theta)\right)dx
\end{eqnarray*}
where
\begin{eqnarray*} 
&&\BARIT b(\BARIT x;\theta)=  -\frac12\PrS\nabla \BARIT U(\BARIT x;\theta)   +  \frac12 \nabla \PrS\PERIOD
\end{eqnarray*}
\end{enumerate}
\end{thm}

\prf
a) Recall the definition of the relative entropy rate $ \mathcal{H}(\PATHP|\tPATHP)$, \cite{KP2013},
$$  \mathcal{H}(\PATHP|\tPATHP) = \frac1h\int\int \mu(x) p^h(x,x')\log \frac{p^h(x,x')}{q^h_\theta(x,x')}dx'dx\PERIOD$$
From the definition of $ p^h, q^h_{\theta}$ and the factorization of $p^h$ given in \VIZ{pFactorize}, we get that
\begin{eqnarray*}
\log\frac{p^h(x,x')}{q^h_\theta(x,x')} &=&  \log\frac{p_1(x,\BARIT x') p_2(x,\hat x'|\BARIT x')}{p_\theta(\BARIT x, \BARIT x') \nu(x'|\BARIT  x')} =  \\
&=& \log\frac{ p_2(x,\hat x'|\BARIT x')}{\nu(x'|\BARIT  x')}   - \frac{1}{2h} \left[ (\BARIT x' -  \COP \Delta x)^{tr} \Sigma_{11}^{-1}(x)(\BARIT x' -   \COP \Delta x)\right] \\
&&+ \log \frac{\BARIT Z^h(\BARIT x;\theta)}{Z_1^h(x)}  
  + \frac{1}{2h} \left[ (\BARIT x' -   \Delta_\theta \BARIT x)^{tr}\BARIT\Sigma^{-1}(\BARIT x;\theta)(\BARIT x' -   \Delta_\theta \BARIT x)\right] \PERIOD
\end{eqnarray*}
Note that the first two terms are $\theta $ independent, therefore they will not contribute to the minimization problem 
$ \nabla_\theta \mathcal H (P|P^\theta) = 0\PERIOD $
Thus 
\begin{eqnarray*}
h\nabla_\theta \mathcal{H}(\PATHP|\tPATHP)  = \nabla_\theta&&\left [ \int\int  \int \mu(x) p_1(x,\BARIT x') p_2(x,\hat x'|\BARIT x')  \right.\times \\
&&\left. \times\left( \log \frac{\BARIT Z^h(\BARIT x;\theta)}{Z_1^h(x)}  
 +\frac{1}{2h} \left[ (\BARIT x' -   \Delta_\theta \BARIT x)^{tr}\BARIT\Sigma^{-1}(\BARIT x;\theta)(\BARIT x' -   \Delta_\theta \BARIT x)\right]\right) d\hat x' d\BARIT x' dx\right] \\
&&= \nabla_\theta \TI + \nabla_\theta \TII\PERIOD
\end{eqnarray*}
\begin{eqnarray*}
\TI &=& \int \int \mu(x) p_1(x,\BARIT x')  \log \frac{\BARIT Z^h(\BARIT x;\theta)}{Z_1^h(x)} 
d\BARIT x' dx =
\frac12 \int  \mu(x)  \log |\BARIT \Sigma(\BARIT x;\theta) (\COP\Sigma(x)\COP^{tr})^{-1}| dx\PERIOD
\end{eqnarray*}
 
\begin{eqnarray*}
\TII &=& \frac{1}{2h} \int \int \mu(x) p_1(x,\BARIT x') \left[ (\BARIT x' -   \Delta_\theta \BARIT x)^{tr}\BARIT\Sigma^{-1}(\BARIT x;\theta)(\BARIT x' -   \Delta_\theta \BARIT x)\right]   d\BARIT x' dx\\
&=&  \frac{1}{2h} \int  \mu(x)\left[ ( \COP \Delta x -   \Delta_\theta \BARIT x)^{tr}\BARIT\Sigma^{-1}(\BARIT x;\theta)( \COP \Delta x -   \Delta_\theta \BARIT x)\right]   dx  + \TIII
\end{eqnarray*}
where the first  term on the right hand side in the previous equality is zero 
and
\begin{eqnarray*}
\TIII &=& \frac{1}{2h} \int \int \mu(x)\frac{1}{Z_1^h(x)} e^{-\frac{1}{2h} (\BARIT x' - \COP \Delta x)^{tr} \Sigma_{11}^{-1}(x)   (\BARIT x' - \COP \Delta x)}\\
&&\qquad\qquad\qquad\qquad\times\left[ (\BARIT x' -    \COP \Delta x)^{tr}\BARIT\Sigma^{-1}(\BARIT x;\theta)(\BARIT x' -   \COP \Delta x)\right]   d\BARIT x' dx\PERIOD
\end{eqnarray*}
Since $\BARIT\Sigma$ is positive definite and symmetric we can write $ \BARIT\Sigma^{-1}(\BARIT x;\theta) =  D^{tr} (\BARIT x;\theta)D(\BARIT x;\theta) $, and perform the change of variables $ s = D(\BARIT x' - \COP \Delta x)$ in the  integral appearing in $III$. Then 

\begin{eqnarray*}
\TIII &=& \frac{1}{2h} \int \int \mu(x)\frac{1}{Z_1^h(x) | D(\BARIT x;\theta)|} e^{-\frac{1}{2h}s^{tr} D^{-tr}(\BARIT x;\theta)\Sigma_{11}^{-1}(x)  D^{-1}(\BARIT x;\theta) s}  s^{tr}s   ds dx\\
&=& \frac{1}{2} \int \mu(x) \mathrm{Tr}\left( \BARIT \Sigma^{-1}(\BARIT x;\theta)\COP\Sigma(x)\COP^{tr} \right) dx\PERIOD
\end{eqnarray*}
We use $\mathrm{Tr}(\cdot)$  to denote  the trace of  a matrix.
Concluding, a local minimizer $\theta^*$ of $ \mathcal{H}(P|P^\theta)$ is given as the solution of 
$\nabla_\theta \mathcal{H}(\PATHP|\tPATHP)=0$ where
\begin{eqnarray*}\label{min}
h\nabla_\theta \mathcal{H}(\PATHP|\tPATHP) = 
\nabla_\theta&&\left[ \frac12 \int  \mu(x)  \log |\BARIT \Sigma(\BARIT x;\theta) (\COP\Sigma(x)\COP^{tr})^{-1}| dx
  \right.+\nonumber \\
&&+\frac{1}{2}\int \mu(x) \mathrm{Tr}\left( \COP \Sigma(x)\COP^{tr} \BARIT\Sigma^{-1}(\BARIT x;\theta)\right) dx +\\
&&+\left.\frac{1}{2h}\int\mu(x) \left(\COP \Delta x- \Delta_{\theta}\BARIT x\right)^{tr} \BARIT\Sigma^{-1}(\BARIT x;\theta)\left(\COP \Delta x- \Delta_{\theta}\BARIT x\right)dx\right]\nonumber\\
&&=  \nabla_\theta\left[ \frac1hA(\theta) +   B(\theta)\right]
\end{eqnarray*}

\medskip
b)
 We  have that  
$ \lim\limits_{h\to 0} \min\limits_{\theta \in \Theta} \left[ \frac1h A(\theta)+B(\theta)\right] = 
\min\limits_{\{\theta \in \Theta: A(\theta)=\frac{m}{2}\}} \left[B(\theta)\right]
$
where  $\frac{m}{2} = \min_\theta A(\theta)   $. This relation can be easily proved  using  optimality arguments as $h \to 0$.
Note that
$A(\theta)= \frac{m}{2} $ when
$$ \BARIT\Sigma(\BARIT x;\theta) = \COP \Sigma(x)\COP^{tr}\PERIOD  $$
The proof is thus completed if  we substitute the above relation in the form of $B(\theta)$.
 \hfill $\square$
%
%
\subsection{ The BBK scheme for Langevin dynamics }
The explicit Euler-Maruyama-Verlet followed by implicit Euler-Maruyama scheme,
also known as BBK scheme, is applied for the discretization of the Langevin dynamics \VIZ{Lang} 
  \begin{equation*}\label{BBK}
 \begin{cases}
  p_{i+1/2} = p_i - F(q_i) \frac{h}{2} - \gamma \mass^{-1}p_i\frac{h}{2} +\sigma \Delta W_i\COMMA\\
  q_{i+1} =q_i +\mass^{-1} p_{i+1/2}h\ \COMMA\\ 
  p_{i+1} = p_{i+1/2} - F(q_{i+1})\frac{h}{2} -\gamma \mass^{-1}p_{i+1}\frac{h}{2} +\sigma \Delta W_{i+1/2}   \end{cases}
   \end{equation*}
 with $\Delta W_i$, $\Delta W_{i+1/2}$  a sequence of i.i.d. Gaussian random vectors with mean zero and covariance $\frac{h}{2} I_{n}$, where for notation simplicity we set $n=3N$. For simplicity we set $\gamma = \gamma I_{n} $ be constant and similarly for $ \sigma = \sigma I_{n} $. 
 Also we set the same mass $M=1$ for all particles.
The discretized Langevin process $(q_i,p_i)$ is a Markov chain, with transition probability from the state $(q,p)$ to $(q',p')$ given by 
\begin{eqnarray*}
P^h(q,p,q',p') = P^h(q'|q,p)P^h(p'|q',q,p)\COMMA
\end{eqnarray*}
where 
\begin{eqnarray*}
&&P^h(q'|q,p) = \frac{1}{Z_0^h}\exp\left\{-\frac{1}{\sigma^2h^3}|q'-\Delta^h_0(q,p)|^2\right\}\;\;\;\;\;  \mbox{and}\\
&&P^h(p'|q',q,p) =  \frac{1}{Z_1^h}\exp\left\{-\frac{1}{\sigma^2h}|p'(1+\frac{\gamma h }{2})-\Delta^h_1(q,q')|^2\right\}\COMMA
\end{eqnarray*}
where we introduce the notation
\begin{eqnarray*}
\Delta^h_0(q,p)  =  q- h\left[p-F(q)\frac{h}{2} + p\frac{h \gamma}{2}\right] \COMMA
\quad \Delta^h_1(q,q') = \frac{1}{h}(q'-q) -\frac{h}{2}F(q')\COMMA
\end{eqnarray*}
and the normalizing constants 
$Z_0^h =\left(  \pi h^3 \sigma^2 \right)^{n/2}$ and $
Z_1^h= \left( \frac{\pi h \sigma^2 }{(1+\gamma h/2)^2}\right)^{n/2}\PERIOD$
The reduced discrete model based on the CG mapping $\COP$  given by the projected  
process $(\COP q_i, \COP p_i) \in \R^{n}\times  \R^{n}$  is approximated by the Markov chain  $ (\BARIT q_i, \BARIT p_i ) \in \R^{m}\times  \R^{m}$ satisfying
\begin{equation*}\label{BBKcg}
\begin{cases} 
 \bp_{i+1/2} = \bp_i - \BF(\bq_i;\theta) \frac{h}{2} - \gamma \bp_i\frac{h}{2} +\sigma \Delta \BARIT{W}_i\COMMA\\
  \bq_{i+1} = \bq_i +  \bp_{i+1/2}h\COMMA\\ 
  \bp_{i+1} =\bp_{i+1/2} - \BF(\bq_{i+1};\theta)\frac{h}{2} - \gamma \bp_{i+1}\frac{h}{2} +\sigma \Delta \BARIT{W}_{i+1/2}
  \end{cases}
\end{equation*}
where $ \Delta \BARIT{W}_i$, $\Delta \BARIT{W}_{i+1/2} \sim \Nn(0,\frac{h}{2}I_{m}) $.
The Markov chain $(\bq_i,\bp_i)$ has transition probabilities
\begin{equation*}
\BP_{\theta}^h(\bq,\bp,\bq',\bp') = \BP_{\theta}^h(\bq'|\bq,\bp)\BP_{\theta}^h(\bp'|\bq',\bq,\bp)\COMMA
\end{equation*}
where 
\begin{eqnarray*}
&&\BP_{\theta}^h(\bq'|\bq,\bp) = \frac{1}{\BARIT Z_0^h}\exp\left\{-\frac{1}{\sigma^2h^3}|\bq'-\BDelta^h_0(\bq,\bp)|^2\right\}\;\;\;\;
 \mbox{ and }\\
&&\BP_{\theta}^h(\bp'|\bq',\bq,\bp) =  \frac{1}{\BARIT Z_1^h}\exp\left\{-\frac{1}{\sigma^2h}|\bp'(1+\frac{\gamma h }{2})-\BDelta^h_1(\bq,\bq')|^2\right\}\COMMA
\end{eqnarray*}
with 
$ \BDelta^h_0(\bq,\bp)=\bq- h\left[\bp-\BF(\bq;\theta)\frac{h}{2} + \bp\frac{h \gamma}{2} \right]$ and 
$ \BDelta^h_1(\bq,\bq') =\frac{1}{h}(\bq'-\bq) -\frac{h}{2}\BF(\bq';\theta)\COMMA$
and normalizing constants 
%
$\BARIT Z_0^h =\left( \pi h^3 \sigma^2 \right)^{m/2}\COMMA$
$\BARIT Z_1^h= \left( \frac{ \pi h \sigma^2 }{(1+\gamma h)^2}\right)^{m/2}\PERIOD
$
We denote with $\nu(q',p'|\bq',\bp',q,p)$ a probability measure associated with the reconstruction map and assume that   is independent of the parameters $\theta$. 
Then the reconstructed chain $(q_i,p_i)$ from $ (\bq_i,\bp_i) $ 
has transition probability
\begin{equation*}
Q_{\theta}^h(q,p, q', p') = \BP_{\theta}^h(\bq,\bp,\bq',\bp')    \nu(q',p'|\bq',\bp',q,p)\PERIOD
\end{equation*}

\begin{thm}\label{prop1Lang}
Given $h>0$ the variational problem for the best-fit of parameters is 
\begin{equation}\label{argminLang}
\mathrm{argmin}_{\theta} \Hh(P^h|Q^h_\theta)  = \mathrm{argmin}_{\theta} \left[ C(\theta)  + D_h(\theta)  \right]\COMMA
\end{equation}
where
 \begin{eqnarray*}
&& C (\theta) = \frac{1}{4}\E_\mu \left[  \left|\sigma^{-1}( \COP F(q)-\BF(\bq;\theta)) \right|^2  \right]\COMMA \\
&&D_h(\theta)= \E_\mu \left[   \frac{1}{Z_0^h}\int  e^{-\frac{1}{\sigma^2h^3} |q'-\Delta^h_0(q,p)|^2 }\left| \sigma^{-1}(\COP F(q')-\BF(\bq';\theta) )\right|^2  dq' \right] \PERIOD
\end{eqnarray*}
Furthermore, as $h\to 0$ we have the optimality condition
\begin{equation}
\lim_{h\to 0} \nabla_{\theta}  \Hh(P^h|Q^h_\theta)   = \frac54 \nabla_{\theta} \E_\mu \left[  \left|\sigma^{-1} (\COP F(q)-\BF(\COP q;\theta) )\right|^2  \right]\PERIOD
\end{equation}
\end{thm}

\begin{prf}
Based on the orthogonality of the coarse graining mapping $\COP$  we can factorize $P^h(q'|q,p)$ and $P^h(p'|q',q,p)$ as
   \begin{equation*}
     P^h(q'|q,p) = \frac{1}{Z_0^h}e^{-\frac{1}{\sigma^2h^3}|\bq'-\COP\Delta^h_0(q,p)|^2}e^{-\frac{1}{\sigma^2h^3}|\tq'-\COP^\perp\Delta^h_0(q,p)|^2}
   \end{equation*}
   and
  \begin{eqnarray*}
    P^h(p'|q',q,p) &\!\!=\!\!&  \frac{1}{Z_1^h}e^{-\frac{1}{\sigma^2h}|\bp'(1+\frac{\gamma h }{2})-\COP\Delta^h_1(q,q')|^2}e^{-\frac{1}{\sigma^2h}|\tp'(1+\frac{\gamma h }{2})-\COP^\perp\Delta^h_1(q,q')|^2}
   \end{eqnarray*}
   where 
  $ \hat {x}$ denotes  $\COP^\perp x = (I-\COP) x$, for $\COP$ a projection,
  $ \COP\Delta^h_0(q,p) = \bq- h\left[\bp- \COP F(q)\frac{h}{2} + \bp\frac{h \gamma}{2} \right] $ and $ \COP\Delta^h_1(q,q')=\frac{1}{h}(\bq'-\bq) -\frac{h}{2}\COP F(q')\COMMA$
    $\COP^\perp\Delta^h_0(q,p) $ and $ \COP^\perp\Delta^h_1(q,q')$ are similarly defined.
From the definition of the RER functional  $ \Hh(P^h|P^h_\theta)$ and the fact that $\nu(q',p'|\bq',\bp',q,p)$ 
is independent of the parameters $\theta$, we have

\begin{eqnarray*}
h\nabla_{\theta} \Hh(P^h|Q^h_\theta) & =   
\nabla_{\theta} \left\{ \idotsint  P^h(q,p,q',p') 
\log\frac{P^h(q,p,q',p')}{P^h_\theta(q,p,q',p')}dq'dp'\mu(dq,dp)\right\}&\\
&= \nabla_{\theta} \TI + \nabla_{\theta} \TII \qquad\qquad\qquad \qquad\qquad\qquad \qquad \qquad\qquad\qquad  &
\end{eqnarray*}
where  
\begin{eqnarray*}
 \TI &=& \frac{1}{\sigma^2h} \frac{ \hat {Z}^h_1}{Z_1^h}\int \dots \int  e^{-\frac{1}{\sigma^2h}|\bp'(1+\frac{\gamma h }{2})-\COP\Delta^h_1(q,q')|^2}P^h(q'|q,p) \\
 &&\qquad\qquad\qquad\qquad\times  \left|\bp'(1+\gamma\frac{h}{2}) - \BDelta^h_1(\bq,\bq')\right|^2  d\bq'd\bp'd\tq' \mu(dq,dp) \COMMA\\
 \TII&=& \frac{1}{\sigma^2h^3}  \frac{ \hat {Z}^h_1}{Z_1^h}\int \dots \int e^{-\frac{1}{\sigma^2h}|\bp'(1+\frac{\gamma h }{2})-\COP\Delta^h_1(q,q')|^2}P^h(q'|q,p) \\
 &&\qquad\qquad\qquad\qquad\times\left|\bq'-\BDelta^h_0(\bq,\bp) \right|^2 d\bq'd\bp'd\tq' \mu(dq,dp)\COMMA
\end{eqnarray*}
and 
$$ 
\hat {Z}^h_1
=\left(\frac{\pi h \sigma^2}{(1+\gamma h/2)^2}\right)^{(n-m)/2} \PERIOD
$$

After integrating with respect to $d\bp'$  and $d\tq' $ we obtain
\begin{eqnarray*}
 \TII 
 &=& \frac{1}{\sigma^2h^3} \frac{1}{Z_0^h} \idotsint e^{-\frac{1}{\sigma^2h^3}
 |\bq'-\COP\Delta^h_0(q,p)|^2 }  \left|\bq'-\BDelta^h_0(\bq,\bp) \right|^2 d\bq'  \mu(dq,dp)\\
 &=& \frac{1}{\sigma^2h^3} \frac{1}{Z_0^h}\idotsint e^{-\frac{1}{\sigma^2h^3} |\bq'-\COP\Delta^h_0(q,p)|^2 }  \left| \COP\Delta^h_0(q,p)-\BDelta^h_0(\bq,\bp) \right|^2 d\bq'  \mu(dq,dp)\\
&& +\; \frac{2}{\sigma^2h^3} \!\!\frac{1}{Z_0^h}\!\!\idotsint \!\!e^{-\frac{1}{\sigma^2h^3}|\bq'-\COP\Delta^h_0(q,p)|^2 }  \!\!\left(\bq'-\BDelta^h_0(\bq,\bp) \right)\!\left( \COP\Delta^h_0(q,p)-\BDelta^h_0(\bq,\bp) \right)\! d\bq'  \mu(dq,dp)\\
&& +\; \frac{1}{\sigma^2h^3} \frac{1}{Z_0^h}\idotsint e^{-\frac{1}{\sigma^2h^3}|\bq'-\COP\Delta^h_0(q,p)|^2 } \! \left|\bq'-\COP\Delta^h_0(q,p) \right|^2 \! d\bq'  \mu(dq,dp)\PERIOD
\end{eqnarray*}
 Therefore
\begin{eqnarray*}
\nabla_{\theta} \TII&=&\nabla_{\theta}\left\{\frac{1}{\sigma^2h^3} \frac{1}{Z_0^h}\idotsint e^{-\frac{1}{\sigma^2h^3} |\bq'-\COP\Delta^h_0(q,p)|^2 }  \left| \COP\Delta^h_0(q,p)-\BDelta^h_0(\bq,\bp) \right|^2 d\bq'  \mu(dq,dp)\right\}\\
&=& \nabla_{\theta}\left\{\frac{1}{\sigma^2h^3}  \iint  \left| \COP\Delta^h_0(q,p)-\BDelta^h_0(\bq,\bp) \right|^2  \mu(dq,dp)\right\} \\
&=& \nabla_{\theta}\left\{\frac{h}{4}  \iint   \left|\sigma^{-1}( \COP F(q)-\BF(\bq;\theta) )\right|^2  \mu(dq,dp)\right\}\PERIOD
\end{eqnarray*}
Similarly,
\begin{eqnarray*}
\nabla_{\theta}  \TI 
&=&\nabla_{\theta}\left\{ h \frac{1}{Z_0^h} \iiint   e^{-\frac{1}{\sigma^2h^3} |q'-\Delta^h_0(q,p)|^2 } 
                         \left| \sigma^{-1}(\COP F(q')-\BF(\bq';\theta) )\right|^2  dq'\mu(dq,dp)\right\}\PERIOD
\end{eqnarray*}
Summarizing all steps we get 
\begin{eqnarray*}
  \nabla_{\theta} \Hh(P^h|Q^h_\theta )& =&  \nabla_{\theta}\left\{ \frac{1}{Z_0^h}\E_\mu\left[\int   e^{\frac{1}{\sigma^2h^3} |q'-\Delta^h_0(q,p)|^2 }\left|\sigma^{-1}( \COP F(q')-\BF(\bq';\theta)) \right|^2  dq'\right] \right\}\\
  &&+ \, \nabla_{\theta}\left\{   \frac{1}{4}  \E_\mu\left[ \left| \sigma^{-1}( \COP F(q)-\BF(\bq;\theta)) \right|^2 \right] \right\}\PERIOD
  \end{eqnarray*}
Furthermore, as $h\to 0^+$ we have that $ \frac{1}{Z_0^h} e^{-\frac{1}{\sigma^2h^3} |\bq'-\COP\Delta^h_0(q,p)|^2 }   \to \delta(q'-q)$ weakly, where $\delta$ denotes the Dirac distribution. Thus 
\begin{eqnarray*}
 \lim_{h\to 0 }\nabla_{\theta} \Hh(P^h|Q^h_\theta ) &=&  \E_\mu\left[\left| \sigma^{-1}(\COP F(q)-\BF(\bq;\theta) )\right|^2  \right] 
                   +   \frac{1}{4}    \E_\mu\left[  \left| \sigma^{-1}(\COP F(q)-\BF(\bq;\theta)) \right|^2 \right]\\
            &=&\frac{5}{4}  \E_\mu\left[  \left| \sigma^{-1}(\COP F(q)-\BF(\COP q;\theta)) \right|^2 \right] \COMMA 
\end{eqnarray*}
which concludes the proof. \hfill $\square$
\end{prf}



 \end{document}